\documentclass[12pt,reqno]{amsart}

\usepackage{amssymb,amsmath,graphicx,amsfonts,euscript}
\usepackage{color}

\newtheorem{thm}{Theorem}[section]

\newtheorem{prop}[thm]{Proposition}
\newtheorem{define}[thm]{Definition}
\newtheorem{rem}[thm]{Remark}

\newtheorem{lemma}[thm]{Lemma}

\newcommand\R{\mathbb {R}}

\numberwithin{equation}{section}
\begin{document}
\title[Ill-posedness for the inhomogeneous Navier-Stokes equations]{Ill-posedness for the 3D inhomogeneous  Navier-Stokes equations in the critical Besov space near $L^6$ framework}

\author[  R. Wan]{ Renhui Wan}
\address{ School of Mathematical Sciences, Nanjing Normal University, Nanjing 210023, China}

\email{rhwanmath@163.com}

\vskip .2in
\begin{abstract}
We prove the ill-posedness for the  3D incompressible inhomogeneous Navier-stokes equations  in critical Besov space. In particular,
a norm inflation happens in finite time with the initial data satisfying
$$\|a_0\|_{\dot{B}_{p,1}^\frac{3}{p}}+\|u_0\|_{\dot{B}_{6,1}^{-\frac{1}{2}}}\le \delta,\ p>6$$
or
$$\|a_0\|_{\dot{B}_{6,1}^\frac{1}{2}}+\|u_0\|_{\dot{B}_{p,1}^{\frac{3}{p}-1}}\le \delta,\ p>6.$$
To obtain the norm inflation,
we   construct a special class of initial data and  introduce a modified pressure. Comparing with the classical Navier-Stokes equations in $L^\infty$ framework, we can obtain the ill-posedness for the inhomogeneous case in near $L^6$ framework.

\end{abstract}

\maketitle

\vskip .2in
\section{Introduction}
\label{Introduction}
In this paper, we consider  the cauchy problem for the 3D incompressible inhomogeneous  Navier-Stokes equations:
\begin{equation} \label{INS1}
\left\{
\begin{aligned}
& \partial_t \rho + {\rm div}(\rho u)= 0,  \\
& \rho\partial_t  u + \rho u\cdot\nabla u-\mu \Delta u+\nabla P =0,\\
& {\rm div}u=0,\\
& (\rho(0,x), u(0,x))=(\rho_0(x), u_0(x)),
\end{aligned}
\right.
\end{equation}
where $(t,x)\in \R^{+}\times \R^3$,  $\rho\in\R$, $u=(u^1,u^2,u^3)\in\R^3$  stand for the density and the velocity field, respectively, $P$ represents the scalar pressure. The constant  $\mu>0$
 is  viscosity  coefficient. $\rho_0 $ and $u_0$ are the initial data satisfying ${\rm div}u_0=0$.
 It is easy to check that the solution $(\rho,u)$ of (\ref{INS1}) is scaling invariance under
\begin{equation}\label{Trans}
(\rho_\lambda,\ u_\lambda)=(\rho(\lambda^2t,\lambda x),\ \lambda u(\lambda^2t,\lambda x)).
\end{equation}
We say a function space is critical means the corresponding norm is invariant under (\ref{Trans}).
\vskip .1in
Lions \cite{Lions} showed (\ref{INS1}) has a global weak solution $(\rho,u)$ with the following initial conditions:
$$0\le \rho_0\in L^\infty,\ \sqrt{\rho}u_0\in L^2.$$
Then Ladyzenskaja and  Solonnikov \cite{LS} obtained the local well-posedness for (\ref{INS1}) with regular data. For the more results on the classical
solution, one can see \cite{Antontsev,Desjardins,Lions} and references therein.
\vskip .1in
Recently, many mathematicians have studied  the well-posedness for the system (\ref{INS1}) in the critical Besov space. Local well-posedness and small data global existence were obtained by Abidi \cite{Ab07} and Danchin \cite{Danchin03}, that is,  \\
local well-posedness:
\begin{equation}\label{well}
\ \rho_0-1\in \dot{B}_{p,1}^\frac{3}{p},\ u_0\in \dot{B}_{p,1}^{\frac{3}{p}-1},\ p<6;
\end{equation}
global well-posedness:
\begin{equation}\label{gwell}
\|\rho_0-1\|_{\dot{B}_{p,1}^\frac{3}{p}}+\|u_0\|_{\dot{B}_{p,1}^{\frac{3}{p}-1}}<\epsilon.
\end{equation}
Particularly, \cite{Ab07,Danchin03}  required the small condition of the initial density and the restriction of $p\in [1,3]$ for the uniqueness, which were removed by \cite{Ab12} and \cite{Danchin12}, respectively.  In fact, all the above results are obtained around the premise that the initial density is near constant 1. It means that
 no vacuum is allowed. Let us  introduce the new unknown $a:=\frac{1}{\rho}-1$,  then (\ref{INS1}) can be rewritten  as follows:
\begin{equation} \label{INS2}
\left\{
\begin{aligned}
& \partial_t a + u\cdot\nabla a= 0,  \\
& \partial_t  u -\mu \Delta u =-u\cdot\nabla u-(1+a)\nabla P+\mu a\Delta u,\\
& {\rm div}u=0,\\
& (a(0,x), u(0,x))=(a_0(x), u_0(x)).
\end{aligned}
\right.
\end{equation}
Paicu and Zhang \cite{PZ} proved the global well-posedness for (\ref{INS2}) with large vertical velocity component (i.e., $u^3$)  with the initial data $(a_0,u_0)\in \dot{B}_{q,1}^\frac{3}{q}\times \dot{B}_{p,1}^{\frac{3}{p}-1}$ satisfying some restrictions on $(p,q)$, which was later improved by the authors in \cite{ZZ}. We refer to \cite{Ab072,Ab11,CZZ,Danchin13,Germain08,HPZ,PZZ} for some other related results. Let us point out that it is not a trivial procedure  to  extend these results  to $L^p$ $(p\ge 6)$ framework, since  there is no effective tool to deal with the nonlinear term $\mu a\Delta u$.
\vskip .1in
When $\rho=$ constant,  (\ref{INS1}) reduces to the classical Navier-Stokes equations. Cannone \cite{Can} and Planchon \cite{Pl} proved global solutions for small data in $\dot{B}_{p,q}^{\frac{3}{p}-1}$ $(p<\infty, q\le \infty)$.  Bourgain and Pavlovic \cite{BP} obtained the ill-posedness in
$\dot{B}_{\infty,\infty}^{-1}$ by proving the solution map is discontinuous in $\dot{B}_{\infty,\infty}^{-1}$. And Germain \cite{Germain082} showed the solution map is not $C^2$ in $\dot{B}_{\infty,q}^{-1}$ $(q>2)$. Yoneda \cite{Y} showed the solution map is not continuous in $\dot{B}_{\infty,q}^{-1}$ $(q>2)$. Very recently, Wang \cite{WangB} obtained the a new ill-posedness in $\dot{B}_{\infty,q}^{-1}$ $(1\le q<2)$. We refer to \cite{100} for the ill-posedness in some Triebel-Lizorkin space and \cite{Cui15,IT14} for other spaces.  We point out that the norm inflation comes from the analysis of  nonlinear term $u\cdot\nabla u$.
\vskip .1in
Roughly speaking,  (\ref{INS2}) is locally well-posedness for the initial data $(a_0,u_0)\in \dot{B}_{p,1}^\frac{3}{p}\times \dot{B}_{p,1}^{\frac{3}{p}-1}$, $p<6$. So a nature question is whether (\ref{INS2}) is well-posedness in the critical Besov space with $p\ge 6$. To the best of our knowledge, similar question has been proposed for the compressible Navier-Stokes, see \cite{CMZ3} for the details. Indeed, the authors \cite{CMZ3}
gave a negative answer to this question, that is,
 the solution   of  the compressible Navier-Stokes equations is ill-posedness  with $p>6$. Very recently, Chen and Wan \cite{CW} proved ill-posedness with the initial  velocity in $L^6$ framework
by using  a new approach to get a norm inflation which depends on  a  decomposition of the density.  Motivated by the above analysis, we will show (\ref{INS2}) is ill-posedness in the critical Besov space. Our main results read:
\begin{thm}\label{t1}
Let $p\in (6,\infty]$. For any $\delta>0$, there exists initial data $(a_0,u_0)$ satisfying
$$\|a_0\|_{\dot{B}_{p,1}^\frac{3}{p}}\le \delta,\ \ \|u_0\|_{\dot{B}_{6,1}^{-\frac{1}{2}}}\le \delta$$
such that a solution $(a,u)$ to the system (\ref{INS2})  satisfies
$$\|u(t)\|_{\dot{B}_{6,1}^{-\frac{1}{2}}}\ge \frac{1}{\delta^\alpha}$$
for some $0<t<\delta$ and $\alpha>0$.
\end{thm}
\begin{thm}\label{t2}
Let  $p\in(6,\infty)$. For any $\delta>0$, there exists initial data $(a_0,u_0)$ satisfying
$$\|a_0\|_{\dot{B}_{6,1}^\frac{1}{2}}\le \delta,\ \ \|u_0\|_{\dot{B}_{p,1}^{\frac{3}{p}-1}}\le \delta$$
such that a solution $(a,u)$ to the system (\ref{INS2}) satisfies
$$\|u(t)\|_{\dot{B}_{p,1}^{\frac{3}{p}-1}}\ge \frac{1}{\delta^\alpha}$$
for some $0<t<\delta$ and $\alpha>0$.
\end{thm}
\begin{rem}\label{r10000}
The idea to the proof of Theorem \ref{t2} can not be applied directly to the case $p=\infty$.  For this case, we
 will give some comments on the barrier and provide a brief framework of the proof  in the Appendix.
\end{rem}
\begin{rem}\label{r1}
The norm inflation for the  classical Navier-stokes equations coming from the nonlinear term $u\cdot\nabla u$ is in $L^\infty$ framework,
while the norm inflation   for
(\ref{INS1}) is in smaller space due to the appearance of  $\mu a\Delta u$.  But we have no idea to extend the  main results to $L^6$ framework.
\end{rem}
\begin{rem}\label{r2000}
In the proof, we only give a priori estimate, which is the key part. We give the structure of the existence in the Appendix. 
\end{rem}
Now, we give the idea of the proof and make some comments on the technics. Firstly, we present our idea.
Like \cite{CW},
the proof is based on a composition of the velocity and a new decomposition of the density (see Section \ref{s3.1}), that is
$$u=U_0+U_1+U_2,$$
$$a=a_0+a_1.$$
Then we  obtain   a norm  inflation  coming from the coupling term $\mu a\Delta u$  yielding a norm inflation of $U_1$, while the corresponding norms of
 $U_0$ and $U_2$ are small.  Secondly, let us show the technics.
\begin{enumerate}
\item[1)] We apply a small trick that a special class of initial velocity is constructed to obtain a large lower bound of the associated norm of $U_1$, see Remark \ref{r2}.
\item[2)] Although we own  the decomposition of the density, we will face the main difficulty coming from the estimate of gradient pressure $(\nabla P)$. As a matter of fact, it seems hard to bound the nonlinear term $\mu |D|^{-2}\nabla {\rm div}(a\Delta u)$. To overcome this difficulty, we introduce a modified pressure $\Pi$ satisfying
$$\nabla \Pi:=\cdot\cdot\cdot-\mu |D|^{-2}\nabla {\rm div}\{a|D|^{-2}\nabla {\rm div}(a\Delta u)\}+\cdot\cdot\cdot.$$
Thanks to that this term   $-\mu |D|^{-2}\nabla {\rm div}\{a|D|^{-2}\nabla {\rm div}(a\Delta u)\}$ admits a good estimate, we can achieve this goal.
\end{enumerate}
\vskip .2in
This paper is organized as follows:\\
In Section \ref{sec:Preliminaries}, we provide some lemmas and the definitions of some spaces. In Section \ref{main1}, we prove Theorem \ref{t1}, while Theorem \ref{t2} is proved in Section \ref{main2}. We split each Section into several steps. In the Appendix, we will consider the case $p=\infty$ in Theorem \ref{t2}.

\vskip .2in
Let us complete this section by describing the notations we shall use in this paper.\\
{\bf Notations} In some places of this paper,  we may use $L^p$ and $\dot{B}_{p,r}^s$  to stand for  $L^p(\R^3)$ and $\dot{B}_{p,r}^s(\R^3)$, respectively.
The uniform constant  $C$, which may be different on different lines, while the constant $C(\cdot)$ means a constant depends on the element(s) in bracket.
$a\approx b$ means $N_1^{-1}a\le b\le N_1a$ for some constant $N_1$, and $a\gg b$ $(a\ll b)$ stands for $a\ge Nb$ $(a\le N^{-1}b)$, where $N$ is a large enough constant.

\section{Preliminaries}
\label{sec:Preliminaries}
In this section, we give some necessary definitions,  propositions and lemmas.
\vskip .1in
The fractional Laplacian operator $|D|^\alpha=(-\Delta)^\frac{\alpha}{2}$ is defined through the Fourier transform, namely,
$$\widehat{|D|^\alpha f}(\xi):=|\xi|^\alpha \widehat{f}(\xi),$$
where the Fourier transform  is given by
$$\widehat{f}(\xi):=\int_{\mathbb{R}^3}e^{-ix\cdot\xi}f(x)dx,\ {\rm or}\ \mathcal{F}(f)(\xi):=\int_{\mathbb{R}^3}e^{-ix\cdot\xi}f(x)dx.$$
Let $\mathfrak{C}=\{\xi\in\mathbb{R}^3,\ \frac{3}{4}\le|\xi|\le\frac{8}{3}\}$. Choose a nonnegative smooth radial function $\varphi$ supported in   $\mathfrak{C}$ such that
$$\sum_{j\in\mathbb{Z}}\varphi(2^{-j}\xi)=1,\ \ \xi\in\mathbb{R}^3\setminus\{0\}.$$
We denote $\varphi_{j}=\varphi(2^{-j}\xi),$ $h=\mathfrak{F}^{-1}\varphi$,  where $\mathfrak{F}^{-1}$ stands for the inverse Fourier transform. Then the dyadic blocks
$\Delta_{j}$ and $S_{j}$ can be defined as follows
$$\Delta_{j}f=\varphi(2^{-j}D)f=2^{3j}\int_{\mathbb{R}^3}h(2^jy)f(x-y)dy,\ \ S_j f=\sum_{k\le j-1}\Delta_k f.$$
 One  easily verifies that with our choice of $\varphi$
$$\Delta_{j}\Delta_{k}f=0\ {\rm if} \ |j-k|\ge2\ \ {\rm and}\ \  \Delta_{j}(S_{k-1}f\Delta_{k}f)=0\  {\rm if}\  |j-k|\ge5.$$
Let us recall the definitions of the  Besov space and Chemin-Lerner type space \cite{CL}.
\begin{define}\label{HB}
 Let $s\in \mathbb{R}$, $(p,q)\in[1,\infty]^2,$ the homogeneous Besov space $\dot{B}_{p,q}^s(\R^3)$ is defined by
$$\dot{B}_{p,q}^{s}(\R^3)=\{f\in \mathfrak{S}'(\R^3);\ \|f\|_{\dot{B}_{p,q}^{s}(\R^3)}<\infty\},$$
where
\begin{equation*}
\|f\|_{\dot{B}_{p,q}^s(\R^3)}=\left\{\begin{aligned}
&\displaystyle (\sum_{j\in \mathbb{Z}}2^{sqj}\|\Delta_{j}f\|_{L^p(\R^3)}^{q})^\frac{1}{q},\ \ \ \ {\rm for} \ \ 1\le q<\infty,\\
&\displaystyle \sup_{j\in\mathbb{Z}}2^{sj}\|\Delta_{j}f\|_{L^p(\R^3)},\ \ \ \ \ \ \ \ {\rm for}\ \ q=\infty,\\
\end{aligned}
\right.
\end{equation*}
and $\mathfrak{S}'(\R^3)$ denotes the dual space of $\mathfrak{S}(\R^3)=\{f\in\mathcal{S}(\mathbb{R}^3);\ \partial^{\alpha}\hat{f}(0)=0;\ \forall\ \alpha\in \ \mathbb{N}^3 $\ {\rm multi-index}\} and can be identified by the quotient space of $\mathcal{S'}/\mathcal{P}$ with the polynomials space $\mathcal{P}$.
\end{define}
\begin{define}\label{CL}
 Let $s\in \mathbb{R}$, $(p,q,r)\in[1,\infty]^3$, $0<T\le \infty$.  The Chemin-Lerner type space $\tilde{L}^r_T\dot{B}_{p,q}^s(\R^3)$ is defined by
$$\tilde{L}^r_T\dot{B}_{p,q}^s(\R^3)=\{f\in \mathfrak{S}'(\R^3);\ \|f\|_{\tilde{L}^r_T\dot{B}_{p,q}^s(\R^3)}<\infty\},$$
where
\begin{equation*}
\|f\|_{\tilde{L}^r_T\dot{B}_{p,q}^s(\R^3)}=\left\{\begin{aligned}
&\displaystyle (\sum_{j\in \mathbb{Z}}2^{sqj}\|\Delta_{j}f\|_{L^r_TL^p(\R^3)}^{q})^\frac{1}{q},\ \ \ \ {\rm for} \ \ 1\le q<\infty,\\
&\displaystyle \sup_{j\in\mathbb{Z}}2^{sj}\|\Delta_{j}f\|_{L^r_TL^p(\R^3)},\ \ \ \ \ \ \ \ {\rm for}\ \ q=\infty.\\
\end{aligned}
\right.
\end{equation*}
It is clear that $\tilde{L}^r_T\dot{B}_{p,r}^s=L^r_T\dot{B}_{p,r}^s$.
\end{define}
Let us  introduce the homogeneous  Bony's  decomposition.
$$uv=T_uv+T_vu+R(u,v),$$
where
$$T_uv=\sum_{j\in \mathbb{Z}}S_{j-1}u\Delta_{j}v,\ \ T_vu=\sum_{j\in \mathbb{Z}}\Delta_{j}uS_{j-1}v,\ \
R(u,v)=\sum_{j\in\mathbb{Z}}\Delta_{j}u\tilde{\Delta}_{j}v,$$
here $\tilde{\Delta}_{j}=\Delta_{j-1}+\Delta_{j}+\Delta_{j+1}.$
\vskip .1in
The following proposition
provides Bernstein type inequalities.
\begin{prop}\label{p1.1}
Let $1\le p\le q\le \infty$. Then for any $\beta,\gamma\in \mathbb{N}^3$, there exists a constant $C$ independent of $f,j$ such that
\begin{enumerate}
\item[1)] If $f$ satisfies
$$
\mbox{supp}\, \widehat{f} \subset \{\xi\in \mathbb{R}^3: \,\, |\xi|
\le \mathcal{K} 2^j \},
$$
 then
$$
\|\partial^\gamma f\|_{L^q(\mathbb{R}^3)} \le C 2^{j|\gamma|  +
j (\frac{3}{p}-\frac{3}{q})} \|f\|_{L^p(\mathbb{R}^3)}.
$$
\item[2)] If $f$ satisfies
\begin{equation*}\label{spp}
\mbox{supp}\, \widehat{f} \subset \{\xi\in \mathbb{R}^3: \,\, \mathcal{K}_12^j
\le |\xi| \le \mathcal{K}_2 2^j \}
\end{equation*}
 then
$$
 \|f\|_{L^p(\mathbb{R}^3)} \le C2^{-j|\gamma|}\sup_{|\beta|=|\gamma|} \|\partial^\beta f\|_{L^p(\mathbb{R}^3)}.
$$
\end{enumerate}
\end{prop}
The standard estimates of the heat equation and transport equation read in the following:
\begin{prop}\label{p1.2}
Let $T>0$, $s\in\R$ and $1\le r\le \infty$. Assume that $u_0\in \dot{B}_{r,1}^s$ and $f\in \tilde{L}^\rho_T\dot{B}_{r,1}^{s-2+\frac{2}{\rho}}$. If $u$
is the solution of the heat equation
\begin{equation*}
\left\{
\begin{aligned}
& \partial_t u -\mu \Delta u= f,  \\
& u(0,x)=u_0(x),
\end{aligned}
\right.
\end{equation*}
with $\mu>0$, then $\forall\ \rho_1\in [\rho,\infty]$, we have
\begin{equation*}
\mu^\frac{1}{\rho_1}\|u\|_{\tilde{L}^{\rho_1}_T\dot{B}_{r,1}^{s+\frac{2}{\rho_1}}}\le C(\|u_0\|_{\dot{B}_{r,1}^s}+\|f\|_{\tilde{L}^\rho_T\dot{B}_{r,1}^{s-2+\frac{2}{\rho}}}).
\end{equation*}
\end{prop}
\begin{prop}\label{p1.3}
Let $T>0$, $s\in (-3\min(\frac{1}{r},\frac{r-1}{r}),1+\frac{3}{r}]$, and $1\le r\le \infty$. Assume that $u$ is the solution of
\begin{equation*}
\left\{
\begin{aligned}
& \partial_t u + v\cdot\nabla u= f,  \\
& u(0,x)=u_0(x),
\end{aligned}
\right.
\end{equation*}
then we have $\forall\ t\in [0,T]$,
$$\|u\|_{\tilde{L}^\infty_T\dot{B}_{r,1}^s}\le (\|u_0\|_{\dot{B}_{r,1}^s}+\|f\|_{L^1_t\dot{B}_{r,1}^s})\exp\{\|\nabla v\|_{L^1_t\dot{B}_{r,1}^\frac{3}{r}}\}.$$
\end{prop}
The Kato-Ponce estimate and some product estimates can be given by
\begin{lemma}\cite{KP}
Let $s>0$, $1\le p,r\le \infty,$ then
\begin{equation}\label{kp}
\|fg\|_{\dot{B}_{p,r}^{s}}\le C\left\{\|f\|_{L^{p_{1}}}\|g\|_{\dot{B}_{p_{2},r}^{s}}+\|g\|_{L^{r_{1}}}\|g\|_{\dot{B}_{r_{2},r}^{s}}\right\},
\end{equation}
where $1\le p_{1},r_{1}\le \infty$ such that $\frac{1}{p}=\frac{1}{p_{1}}+\frac{1}{p_{2}}=\frac{1}{r_{1}}+\frac{1}{r_{2}}$.
\end{lemma}
\begin{lemma} \cite{CW} Let
$$1\le s,s_1,s_2,s_{i1},s_{i2}\le \infty,\ 3<r<\infty,\ 3< q<6,\ |\alpha|\ge 1,$$
$$\frac{1}{s}=\frac{1}{s_1}+\frac{1}{s_2}=\frac{1}{s_{i1}}+\frac{1}{s_{i2}},\ i=1,2,$$
$$\frac{3}{p_0}+\frac{3}{r}>1,\ p_0>6.$$
Then   the following estimates  hold:\\
($a$) \begin{equation}\label{p1}
\|fg\|_{\tilde{L}^s_T\dot{B}_{r,1}^\frac{3}{r}}\le C\left(\|g\|_{\tilde{L}^{s_{11}}_TL^\infty}\|f\|_{L^{s_{12}}_T\dot{B}_{r,1}^\frac{3}{r}}
+\|\partial^\alpha g\|_{L^{s_{21}}_TL^\infty}\|f\|_{\tilde{L}^{s_2}_T\dot{B}_{r,1}^{\frac{3}{r}-|\alpha|}}
\right);
\end{equation}
($b$) \begin{equation}\label{p2}
\|fg\|_{\tilde{L}^s_T\dot{B}_{r,1}^{\frac{3}{r}-1}}\le C\left(\|f\|_{\tilde{L}^{s_{11}}_TL^\infty}\|g\|_{L^{s_{12}}_T\dot{B}_{r,1}^{\frac{3}{r}-1}}
+\| f\|_{L^{s_{21}}_T\dot{B}_{p_0,1}^{\frac{3}{p_0}+2}}\|g\|_{\tilde{L}^{s_2}_T\dot{B}_{r,1}^{\frac{3}{r}-3}}
\right);
\end{equation}
($c$) \begin{equation}\label{p3}
\|fg\|_{\tilde{L}^s_T\dot{B}_{r,1}^{\frac{3}{r}-1}}\le C\|f\|_{\tilde{L}^{s_1}_T\dot{B}_{p_0,1}^\frac{3}{p_0}}\|g\|_{\tilde{L}^{s_2}_T\dot{B}_{r,1}^{\frac{3}{r}-1}}.
\end{equation}
\end{lemma}
For the readers' convenience, we refer to \cite{BCD} for more details about the Besov space.
\vskip .4in
\section{Proof of Theorem \ref{t1}}
\label{main1}
\subsection{Reformulation of the equation}
\label{s3.1}
 From (\ref{INS2}), we have
\begin{equation}\label{P}
\nabla P=|D|^{-2}\nabla {\rm div}(a\nabla P)+|D|^{-2}\nabla {\rm div}(u\cdot\nabla u)-\mu |D|^{-2}\nabla {\rm div}(a\Delta u).
\end{equation}
As the previous comments in Section \ref{Introduction}, we require   a modified pressure given by
\begin{equation*}
\Pi:=P+\mu |D|^{-2}{\rm div}(a\Delta u).
\end{equation*}
Then one  gets
\begin{equation}\label{MP}
\nabla \Pi=|D|^{-2}\nabla {\rm div}(a\nabla \Pi)-\mu |D|^{-2}\nabla {\rm div}\left\{a|D|^{-2}\nabla {\rm div}(a\Delta u)\right\}+|D|^{-2}\nabla {\rm div}(u\cdot\nabla u).
\end{equation}
So we can write (\ref{INS2}) as
\begin{equation} \label{INS3}
\left\{
\begin{aligned}
& \partial_t a + u\cdot\nabla a= 0,\   {\rm div}u=0 \\
& \partial_t  u -\mu \Delta u =-u\cdot\nabla u-\nabla \Pi-a\nabla P\\
&\hspace{15mm}+\mu a\Delta u+\mu |D|^{-2}\nabla{\rm div}(a\Delta u),\\
& (a(0,x), u(0,x))=(a_0(x), u_0(x)).
\end{aligned}
\right.
\end{equation}
Applying Duhamel principle to (\ref{INS3})$_2$, we get
\begin{equation*}
\begin{aligned}
u(t,x)=&e^{\mu\Delta t}u_0+\int_0^t e^{\mu \Delta(t-\tau)}\{-u\cdot\nabla u-\nabla \Pi-a\nabla P\\
&+\mu a\Delta u+\mu |D|^{-2}\nabla{\rm div}(a\Delta u)\}d\tau.
\end{aligned}
\end{equation*}
Denote
\begin{equation*}
\begin{aligned}
U_0(t)=&e^{\mu \Delta t}u_0,\\
U_1(t)=&\mu\int_0^t e^{\mu \Delta (t-\tau)}\{a_0\Delta U_0+|D|^{-2}\nabla {\rm div}(a_0\Delta U_0) \}d\tau,\\
U_2(t)=&\int_0^t e^{\mu \Delta (t-\tau)}\{-u\cdot\nabla u-\nabla \Pi-a \nabla P+F_1+F_2\}d\tau,
\end{aligned}
\end{equation*}
where
\begin{equation*}
\begin{aligned}
F_1=&\mu (a_1\Delta u+a_0\Delta(U_1+U_2)),\\
F_2=&\mu|D|^{-2}\nabla {\rm div}\{a_1\Delta u+a_0\Delta (U_1+U_2)\}.
\end{aligned}
\end{equation*}
Now, we can decompose $u(t,x)$ and $a(t,x)$ as
$$u(t,x)=U_0(t,x)+U_1(t,x)+U_2(t,x)$$
and
$$a(t,x)=a_0(x)+a_1(t,x)$$
where $a_1(t)$ satisfies the generalized transport equation given by
\begin{equation} \label{TE}
\left\{
\begin{aligned}
& \partial_t a_1 + u\cdot\nabla (a_0+a_1)= 0, \\
&  a_1(0,x)=0.
\end{aligned}
\right.
\end{equation}

\subsection{The choice of initial data}
\label{s3.2}
Due to supp $\varphi(2^8\xi)\subset \{\frac{3}{4}\times 2^{-8}\le |\xi|\le \frac{8}{3}\times2^{-8}\}$, we get there exists a positive constant $A>0$ such that
at least one of the following two inequalities holds:
\begin{equation}\label{f1}
\varphi(2^8\xi)\frac{\xi_1^2+\xi_3^2}{|\xi|^2}\ge A\varphi(2^8\xi),
\end{equation}
$$\varphi(2^8\xi)\frac{\xi_2^2+\xi_3^2}{|\xi|^2}\ge A\varphi(2^8\xi).$$
Without loss of generality, we assume (\ref{f1}) holds in this present article. Let $C(N)=2^{\frac{N}{2}(\frac{1}{2}-\frac{3}{p}-2\epsilon_1)}$ for some $0<\epsilon_1<\frac{1}{2}(\frac{1}{2}-\frac{3}{p})$ and  $p>6$, where  $N>0$ determined later is a sufficiently large constant leading to $\frac{N}{C(N)}\ll 1$.
We construct the initial data $(a_0,u_0)$ as follows:
\begin{equation*}
\begin{aligned}
\widehat{a_0}(\xi)=&\frac{1}{C(N)}\sum_{k=100}^N2^{-k\frac{3}{p}}\left(\hat{\phi}(\xi-2^ke_1)+\hat{\phi}(\xi+2^ke_1)\right),\\
\widehat{u_0}(\xi)=&\frac{1}{C(N)}\sum_{k=100}^N2^{\frac{k}{2}}\left((\hat{\phi}(\xi+2^ke_1)-\hat{\phi}(\xi-2^ke_1)\right)\frac{1}{|\xi|}
\left(
  \begin{array}{c}
   \xi_2 \\
   -\xi_1 \\
    0 \\
  \end{array}
\right),
\end{aligned}
\end{equation*}
where $e_1=(1,0,0)$ and  $\hat{\phi}$ is a smooth, radial and nonnegative function in $\R^3$ satisfying
\begin{equation*}
\hat{\phi}=\left\{
\begin{aligned}
& 1\ \ {\rm for}\ \ |\xi|\le 1,\\
& 0\ \ {\rm for}\ \ |\xi|\ge 2.\\
\end{aligned}
\right.
\end{equation*}
One can see  $a_0$ is a real valued function, while $u_0$ is a real vector-valued function by observing
$$u_0(x)=\frac{2}{C(N)}\sum_{k=100}^N2^{\frac{k}{2}}\left\{-\mathcal{R}_2(\phi(x)\sin(2^kx_1)),\ \mathcal{R}_1(\phi(x)\sin(2^kx_1)),\ 0\right\},$$
where $\mathcal{R}_i$ is the Riesz transform defined by
$$\widehat{\mathcal{R}_if}(\xi):=\frac{-i\xi_i}{|\xi|}\widehat{f}(\xi).$$
  One can also check the following estimates hold, i.e.,
\begin{equation}\label{3.1}
\|a_0\|_{L^\infty}\le \frac{C}{C(N)},\  \|a_0\|_{\dot{B}_{p,1}^\frac{3}{p}}\le \frac{CN}{C(N)},\ \|u_0\|_{\dot{B}_{6,1}^{-\frac{1}{2}}}\le \frac{CN}{C(N)}
\end{equation}
and
\begin{equation*}
\|a_0\|_{\dot{B}_{r,1}^s}\le \frac{C2^{N(s-\frac{3}{p})}}{C(N)},\ \|u_0\|_{\dot{B}_{r_1,1}^{s_1}}\le \frac{C2^{N(s_1+\frac{1}{2})}}{C(N)},
\end{equation*}
with $(r,r_1)\in [1,\infty]^2$ and $s >\frac{3}{p}$, $s_1>-\frac{1}{2}$.
\vskip .1in
Next, a lemma  is given.
\begin{lemma}\label{l1}
Let $p>6$. Then there exist some positive constants $q$, $p_0$, $\epsilon$, $\epsilon_1$ satisfying
\begin{equation}\label{cond}
\left\{
\begin{aligned}
& 1<\frac{3}{q}+\frac{3}{p_0}<\frac{5}{4}+\frac{3}{2p}-\frac{3}{p_0}+2\epsilon-3\epsilon_1,\\
& \frac{3}{p_0}<\frac{1}{4}+\frac{3}{2p}-\epsilon_1,\ p_0\in (6,p),\\
& \frac{3}{q}<1+3\epsilon-2\epsilon_1,\ q\in (3,6),\\
& 0<2\epsilon<2\epsilon_1<\frac{1}{2}-\frac{3}{p}.
\end{aligned}
\right.
\end{equation}
\end{lemma}
\begin{rem}\label{rr1}
We give the following explanations of the limitations in (\ref{cond}). Let  $T=2^{-2(1+\epsilon)N}$, then we have
\begin{equation}\label{cond11}
\left\{
\begin{aligned}
& \frac{2^{N(\frac{3}{p_0}-\frac{3}{p})}}{C(N)}\ll 1\ \ \Longleftrightarrow\ \ \frac{3}{p_0}<\frac{1}{4}+\frac{3}{2p}-\epsilon_1,\\
& \frac{T^\frac{3}{2}2^{N(\frac{3}{q}-\frac{3}{p}+\frac{5}{2})}}{C(N)^2}\ll1\ \ \Longleftrightarrow\ \ \frac{3}{q}<1+3\epsilon-2\epsilon_1,\\
& \frac{T2^{N(\frac{3}{q}
+\frac{6}{p_0}-\frac{6}{p}+\frac{3}{2})}}{C(N)^3}\ll1\ \ \Longleftrightarrow\ \ \frac{3}{q}+\frac{3}{p_0}<\frac{5}{4}+\frac{3}{2p}-\frac{3}{p_0}+2\epsilon-3\epsilon_1.
\end{aligned}
\right.
\end{equation}
Actually,  we  assume the conditions on the right hand side of (\ref{cond11}) to  ensure the conditions on the left hand side which is required in  our proof.
Furthermore, we  use $\frac{3}{q}+\frac{3}{p_0}>1$  to ensure some product estimates like (\ref{p3}).
The choice of $C(N)$ needs $2\epsilon_1<\frac{1}{2}-\frac{3}{p}$, while $\epsilon<\epsilon_1$ ensures the norm inflation of $U_1$ (see the end of subsection \ref{s3.3}).
\end{rem}
Lemma \ref{l1} can be proved easily, here we use the following  example in this article:
$$\epsilon_1=\frac{1}{4}(\frac{1}{2}-\frac{3}{p}),\ \epsilon=\frac{1}{5}(\frac{1}{2}-\frac{3}{p})$$
and
$$\frac{3}{p_0}=\frac{1}{16}+\frac{21}{8p},\ \forall\ \frac{3}{q}\in (\frac{15}{16}-\frac{21}{8p},\ \frac{19}{20}-\frac{27}{10p}).$$
This gives that
\begin{equation}\label{100}
2(\epsilon_1-\epsilon)=\frac{1}{10}(\frac{1}{2}-\frac{3}{p}),\ \frac{1}{2}(\frac{1}{2}-\frac{3}{p}-2\epsilon_1)=\frac{1}{4}(\frac{1}{2}-\frac{3}{p}).
\end{equation}
\subsection{The analysis of $U_1$}
\label{s3.3}
Let $V^j$ be the $j$-th component of the vector $V$. Thanks to $\dot{B}_{6,1}^{-\frac{1}{2}}\hookrightarrow \dot{B}_{\infty,\infty}^{-1}$, we have
$$\|U_1\|_{\dot{B}_{6,1}^{-\frac{1}{2}}}\ge c\|U_1\|_{\dot{B}_{\infty,\infty}^{-1}}\ge  c|\int \varphi(2^8\xi)\widehat{U_1}(\xi)d\xi|
\ge c|\int \varphi(2^8\xi)\widehat{U_1^2}(\xi)d\xi|
.$$
Let us give the second component $U_1^2$ of $U_1$:
\begin{equation*}
\begin{aligned}
U_1^2=&\mu \int_0^t e^{\mu\Delta(t-\tau)}\{a_0\Delta U_0^2+|D|^{-2}\partial_2{\rm div}(a_0\Delta U_0)\}d\tau\\
=& \mu \int_0^t e^{\mu\Delta(t-\tau)}(1+\partial_2^2|D|^{-2})(a_0\Delta U_0^2)d\tau\\
&+\mu \int_0^t e^{\mu\Delta(t-\tau)}|D|^{-2}\partial_1\partial_2(a_0\Delta U_0^1)d\tau\\
=:&U_{11}^2+U_{12}^2.
\end{aligned}
\end{equation*}
So
$$\|U_1\|_{\dot{B}_{6,1}^{-\frac{1}{2}}}\ge |B_1|-|B_2|,$$
where
$$B_1=\int \varphi(2^8\xi)\widehat{U_{11}^2}(\xi)d\xi,\ B_2=\int \varphi(2^8\xi)\widehat{U_{12}^2}(\xi)d\xi.$$
Now, we give the estimates of $B_1$ and $B_2$.
\vskip .1in
$\bullet$ {\bf The estimate of $B_2$}
\vskip .1in
Using some facts of Fourier transform, we have
$$B_2=-\mu \int \varphi(2^8\xi)\int_0^t e^{-\mu(t-\tau)|\xi|^2}\frac{\xi_1\xi_2}{|\xi|^2}\mathcal{F}(a_0\Delta U_0^1)d\tau d\xi.$$
Thanks to the construction of initial data, we get
\begin{equation*}
\begin{aligned}
\mathcal{F}(a_0\Delta U_0^1)=& -\int \widehat{a_0}(\xi-\eta)|\eta|^2\widehat{U_0^1}(\eta)d\eta\\
=& -\frac{C}{C(N)^2}\sum_{k=100}^N2^{k(\frac{1}{2}-\frac{3}{p})}\int e^{-\mu |\eta|^2\tau}\eta_2|\eta|
A(\xi,\eta,k)d\eta,
\end{aligned}
\end{equation*}
where
$$A(\xi,\eta,k):=-\hat{\phi}(\xi-\eta+2^ke_1)\hat{\phi}(\eta-2^ke_1)+\hat{\phi}(\xi-\eta-2^ke_1)\hat{\phi}(\eta+2^ke_1).$$
Thus $B_2$ can be given by
\begin{equation*}
\begin{aligned}
B_2=\frac{C}{C(N)^2}\sum_{k=100}^N2^{k(\frac{1}{2}-\frac{3}{p})}\int\int
\varphi(2^8\xi)\frac{\xi_1\xi_2}{|\xi|^2}\eta_2|\eta|A(\xi,\eta,k)\int_0^t e^{-\mu((t-\tau)|\xi|^2+\tau|\eta|^2)}d\tau d\xi d\eta.
\end{aligned}
\end{equation*}
Due to $|\xi|\thickapprox 1$, $|\eta_2|\thickapprox1$ and $|\eta|\thickapprox 2^k$, we get
\begin{equation}\label{1001}
\begin{aligned}
|B_2|\le &\frac{Ct}{C(N)^2}\sum_{k=100}^N2^{k(\frac{3}{2}-\frac{3}{p})}\int\int
\varphi(2^8\xi)|A(\xi,\eta,k)| d\xi d\eta\\
\le& \frac{Ct}{C(N)^2}\sum_{k=100}^N2^{k(\frac{3}{2}-\frac{3}{p})}\le \frac{Ct}{C(N)^2}2^{N(\frac{3}{2}-\frac{3}{p})}.
\end{aligned}
\end{equation}
\vskip .1in
$\bullet$ {\bf The estimate of $B_1$}
\vskip .1in
We will show the large lower bound of $B_1$ which yields the norm inflation of the solution. One can easily obtain
\begin{equation*}
B_1=\mu \int \varphi(2^8 \xi)\frac{\xi_1^2+\xi_3^2}{|\xi|^2} \int_0^t e^{-\mu|\xi|^2(t-\tau)}\mathcal{F}(a_0\Delta U_0^2)d\tau d\xi.
\end{equation*}
By a similar way as before, we can obtain
\begin{equation*}
\begin{aligned}
\mathcal{F}(a_0\Delta U_0^2)=& -\int \widehat{a_0}(\xi-\eta)|\eta|^2\widehat{U_0^2}(\eta)d\eta\\
=&\frac{C}{C(N)^2}\sum_{k=100}^N2^{k(\frac{1}{2}-\frac{3}{p})}\int A(\xi,\eta,k) \eta_1|\eta|e^{-\mu |\eta|^2\tau}d\eta.
\end{aligned}
\end{equation*}
Hence
$$B_1=-\frac{C}{C(N)^2}\sum_{k=100}^N2^{k(\frac{1}{2}-\frac{3}{p})}
\int\int \varphi (2^8 \xi) \frac{\xi_1^2+\xi_3^2}{|\xi|^2}\eta_1|\eta|A(\xi,\eta,k)\mathfrak{A}(t,\xi,\eta)d\xi d\eta,
$$
where
$$\mathfrak{A}(t,\xi,\eta):=\int_0^t e^{-\mu(|\xi|^2(t-\tau)+|\eta|^2\tau)}d\tau.$$
\begin{rem}\label{r2}
Due to the construction of initial data, we obtain two negative terms:
$$\eta_1A(\xi,\eta,k)=-\eta_1\hat{\phi}(\xi-\eta+2^ke_1)\hat{\phi}(\eta-2^ke_1)+\eta_1\hat{\phi}(\xi-\eta-2^ke_1)\hat{\phi}(\eta+2^ke_1).$$
In the following proof, we only use one of them.
\end{rem}
Applying the Taylor expansion $e^x=\sum_{r\ge 0}\frac{x^r}{r!}$, $|\xi|\thickapprox 1$ and $|\eta|\thickapprox 2^k$,  we get
\begin{equation}\label{3.3}
\mathfrak{A}(t,\xi,\eta)=\frac{e^{-\mu t|\eta|^2}-e^{-\mu t|\xi|^2}}{\mu (|\xi|^2-|\eta|^2)}=t+\mathcal{O}(t^2|\eta|^2)
\end{equation}
when $t2^{2N}<1.$
Thanks to (\ref{3.3}) and (\ref{f1}), one has
\begin{equation}\label{1002}
\begin{aligned}
|B_1|\ge& \frac{C}{C(N)^2}\sum_{k=100}^N2^{k(\frac{1}{2}-\frac{3}{p})}
\int\int \varphi (2^8 \xi) \frac{\xi_1^2+\xi_3^2}{|\xi|^2}\eta_1|\eta|\\
&\times
\hat{\phi}(\xi-\eta+2^ke_1)\hat{\phi}(\eta-2^ke_1)\{t+\mathcal{O}(t^2|\eta|^2)\}d\xi d\eta\\
\ge &  \frac{c}{C(N)^2}\sum_{k=100}^N2^{k(\frac{1}{2}-\frac{3}{p})}(t2^{2k}-\mathcal{O}(t^22^{4k}))\\
\ge& \frac{ct}{C(N)^2}\sum_{k=100}^N2^{k(\frac{5}{2}-\frac{3}{p})}-\frac{Ct^2}{C(N)^2}\sum_{k=100}^N2^{k(\frac{9}{2}-\frac{3}{p})}\\
=& \frac{ct2^{N(\frac{5}{2}-\frac{3}{p})}}{C(N)^2}-\frac{Ct^22^{N(\frac{9}{2}-\frac{3}{p})}}{C(N)^2}.
\end{aligned}
\end{equation}
Choosing $t=T_0:=2^{-2(1+\epsilon)N}$, $0<\epsilon<\epsilon_1$, which ensures $t2^{2N}<1$, and combining   with (\ref{1001}), (\ref{1002}) yields
\begin{equation}\label{LB}
\begin{aligned}
\|U_1(t)\|_{\dot{B}_{6,1}^{-\frac{1}{2}}}\ge& |B_1|-|B_2|\\
\ge &  \frac{ct2^{N(\frac{5}{2}-\frac{3}{p})}}{C(N)^2}-\frac{Ct^22^{N(\frac{9}{2}-\frac{3}{p})}}{C(N)^2}-\frac{Ct2^{N(\frac{3}{2}-\frac{3}{p})}}{C(N)^2}\\
\ge& \frac{ct2^{N(\frac{5}{2}-\frac{3}{p})}}{2C(N)^2}-\frac{Ct^22^{N(\frac{9}{2}-\frac{3}{p})}}{C(N)^2}
\ge \frac{c}{4}2^{2(\epsilon_1-\epsilon)N}.
\end{aligned}
\end{equation}
\subsection{The analysis of $U_2$}
\label{s3.4}
 Let $(p,p_0,q)$ be given as in  Lemma \ref{l1}.  Let $0\le T\le T_0$. We split the analysis into five steps.\\
{\bf Step 1. Some estimates of $U_1$} We provide some estimates of $U_1$  which will be used in the following proof.
$$\|U_1\|_{L^r_T L^\infty}\le T^\frac{1}{r}\int_0^T \|a_0\|_{L^\infty}\||\xi|^2\widehat{U_0}\|_{L^1}d\tau\le\frac{CT^{1+\frac{1}{r}}2^{\frac{5}{2}N}}{C(N)^2},$$
and
\begin{equation*}
\begin{aligned}
&\|U_1\|_{L^1_T\dot{B}_{q_1,1}^{\frac{3}{q_1}+1}}+\|U_1\|_{\tilde{L}^2_T\dot{B}_{q_1,1}^\frac{3}{q_1}}\\
\le& T^\frac{1}{2}(\|U_1\|_{\tilde{L}^2_T\dot{B}_{q_1,1}^{\frac{3}{q_1}+1}}+\|U_1\|_{\tilde{L}^\infty_T\dot{B}_{q_1,1}^\frac{3}{q_1}})\\
\le& CT^\frac{1}{2}\int_0^T \|a_0\Delta U_0\|_{\dot{B}_{q_1,1}^\frac{3}{q_1}}d\tau\\
\le& CT^\frac{3}{2}(\|a_0\|_{L^\infty}\|\Delta U_0\|_{\tilde{L}^\infty_T\dot{B}_{q_1,1}^\frac{3}{q_1}}+\|a_0\|_{\dot{B}_{q_1,1}^\frac{3}{q_1}}\||\xi|^2\widehat{U_0}\|_{L^\infty_TL^1})\\
\le&\frac{CT^\frac{3}{2}2^{N(\frac{3}{q_1}+\frac{5}{2})}}{C(N)^2},
\end{aligned}
\end{equation*}
where $q_1=q\  {\rm or}\  p_0$.\\
{\bf Step 2. The estimate of $\|u\cdot\nabla u\|_{L^1_T\dot{B}_{q,1}^{\frac{3}{q}-1}}$}
Thanks to ${\rm div}u=0$ and Bernstein inequality, it suffices to bound $\|u\otimes u\|_{L^1_T\dot{B}_{q,1}^\frac{3}{q}}$. Using the decomposition $u=U_0+U_1+U_2$, we can split this estimate into six parts. Applying (\ref{kp}), one has
\begin{equation*}
\begin{aligned}
\|U_0\otimes U_0\|_{L^1_T\dot{B}_{q,1}^\frac{3}{q}}\le& C\|U_0\|_{L^\infty_TL^\infty}\|U_0\|_{\tilde{L}^1_T\dot{B}_{q,1}^\frac{3}{q}}
 \le \frac{CT2^{N(\frac{3}{q}+1)}}{C(N)^2},\\
\|U_0\otimes U_1\|_{L^1_T\dot{B}_{q,1}^\frac{3}{q}}\le& CT^\frac{1}{2}(\|U_0\|_{L^\infty_TL^\infty}\|U_1\|_{\tilde{L}^2_T\dot{B}_{q,1}^\frac{3}{q}}
+\|U_1\|_{L^\infty_TL^\infty}\|U_0\|_{\tilde{L}^2_T\dot{B}_{q,1}^\frac{3}{q}})
 \le \frac{CT^22^{N(\frac{3}{q}+3)}}{C(N)^3},\\
\|U_1\otimes U_1\|_{L^1_T\dot{B}_{q,1}^\frac{3}{q}}\le& C\|U_1\|_{L^2_TL^\infty}\|U_1\|_{\tilde{L}^2_T\dot{B}_{q,1}^\frac{3}{q}}\le \frac{CT^32^{N(\frac{3}{q}+5)}}{C(N)^4}.
\end{aligned}
\end{equation*}
Using (\ref{p1}),
\begin{equation*}
\begin{aligned}
\|U_0\otimes U_2\|_{L^1_T\dot{B}_{q,1}^\frac{3}{q}}\le& C(\|U_0\|_{L^2_TL^\infty}\|U_2\|_{\tilde{L}^2_T\dot{B}_{q,1}^\frac{3}{q}}
+\|\nabla U_0\|_{L^1_TL^\infty}\|U_2\|_{\tilde{L}^\infty_T\dot{B}_{q,1}^{\frac{3}{q}-1}}
)\\
\le& \frac{CT^\frac{1}{2}2^\frac{N}{2}}{C(N)}(\|U_2\|_{\tilde{L}^\infty_T\dot{B}_{q,1}^{\frac{3}{q}-1}}+\|U_2\|_{\tilde{L}^2_T\dot{B}_{q,1}^\frac{3}{q}}),\\
\|U_1\otimes U_2\|_{L^1_T\dot{B}_{q,1}^\frac{3}{q}}\le& C(\|U_1\|_{\tilde{L}^2_TL^\infty}\|U_2\|_{\tilde{L}^2_T\dot{B}_{q,1}^\frac{3}{q}}
+\|\nabla U_1\|_{L^1_TL^\infty}\|U_2\|_{\tilde{L}^\infty_T\dot{B}_{q,1}^{\frac{3}{q}-1}}
)\\
\le& \frac{CT^\frac{3}{2}2^\frac{5N}{2}}{C(N)^2}(\|U_2\|_{\tilde{L}^\infty_T\dot{B}_{q,1}^{\frac{3}{q}-1}}+\|U_2\|_{\tilde{L}^2_T\dot{B}_{q,1}^\frac{3}{q}}).
\end{aligned}
\end{equation*}
In reality, we have applied Proposition \ref{p1.2} with $f=0$ to some  estimates of $U_0$.  Using (\ref{kp}) again, with $\dot{B}_{q,1}^\frac{3}{q}\hookrightarrow L^\infty$, we have
$$\|U_2\otimes U_2\|_{L^1_T\dot{B}_{q,1}^\frac{3}{q}}\le C\|U_2\|_{\tilde{L}^2_T\dot{B}_{q,1}^\frac{3}{q}}^2.$$
Thus we get
\begin{equation}\label{aim1}
\|u\cdot\nabla u\|_{L^1_T\dot{B}_{q,1}^{\frac{3}{q}-1}}
\le \frac{CT2^{N(\frac{3}{q}+1)}}{C(N)^2}+\frac{C}{C(N)}
(\|U_2\|_{\tilde{L}^\infty_T\dot{B}_{q,1}^{\frac{3}{q}-1}}+\|U_2\|_{\tilde{L}^2_T\dot{B}_{q,1}^\frac{3}{q}})
+C\|U_2\|_{\tilde{L}^2_T\dot{B}_{q,1}^\frac{3}{q}}^2.
\end{equation}
{\bf Step 3. The estimate of $\|F_1\|_{L^1_T\dot{B}_{q,1}^{\frac{3}{q}-1}}$ and $\|F_2\|_{L^1_T\dot{B}_{q,1}^{\frac{3}{q}-1}}$}
Thanks to  the product estimate (\ref{p3}) with $\frac{3}{q}+\frac{3}{p_0}>1$, we can deduce that
\begin{equation*}
\begin{aligned}
\|a_1\Delta U_0\|_{L^1_T\dot{B}_{q,1}^{\frac{3}{q}-1}}\le & C\|a_1\|_{\tilde{L}^\infty_T\dot{B}_{p_0,1}^\frac{3}{p_0}}\|U_0\|_{L^1_T\dot{B}_{q,1}^{\frac{3}{q}+1}}\le \frac{CT2^{N(\frac{3}{q}+\frac{3}{2})}}{C(N)}
\|a_1\|_{\tilde{L}^\infty_T\dot{B}_{p_0,1}^\frac{3}{p_0}},\\
\|a_1\Delta U_1\|_{L^1_T\dot{B}_{q,1}^{\frac{3}{q}-1}}\le & C\|a_1\|_{\tilde{L}^\infty_T\dot{B}_{p_0,1}^\frac{3}{p_0}}\|U_1\|_{L^1_T\dot{B}_{q,1}^{\frac{3}{q}+1}}\le \frac{CT^\frac{3}{2}2^{N(\frac{3}{q}+\frac{5}{2})}}{C(N)^2}
\|a_1\|_{\tilde{L}^\infty_T\dot{B}_{p_0,1}^\frac{3}{p_0}},\\
\|a_1\Delta U_2\|_{L^1_T\dot{B}_{q,1}^{\frac{3}{q}-1}}\le & C\|a_1\|_{\tilde{L}^\infty_T\dot{B}_{p_0,1}^\frac{3}{p_0}}\|U_2\|_{L^1_T\dot{B}_{q,1}^{\frac{3}{q}+1}},\\
\|a_0\Delta U_1\|_{L^1_T\dot{B}_{q,1}^{\frac{3}{q}-1}}\le &\|a_0\|_{\dot{B}_{p_0,1}^\frac{3}{p_0}}\|U_1\|_{L^1_T\dot{B}_{q,1}^{\frac{3}{q}+1}}
\le\frac{CT^\frac{3}{2}2^{N(\frac{3}{q}+\frac{3}{p_0}-\frac{3}{p}+\frac{5}{2})}}{C(N)^3}.
\end{aligned}
\end{equation*}
Applying (\ref{p2}), we can get
\begin{equation*}
\begin{aligned}
\|a_0\Delta U_2\|_{L^1_T\dot{B}_{q,1}^{\frac{3}{q}-1}}\le&
 C(\|a_0\|_{L^\infty}\|U_2\|_{L^1_T\dot{B}_{q,1}^{\frac{3}{q}+1}}+T\|a_0\|_{\dot{B}_{p_0,1}^{\frac{3}{p_0}+2}}
 \|U_2\|_{\tilde{L}^\infty_T\dot{B}_{q,1}^{\frac{3}{q}-1}})\\
 \le& C\frac{T2^{N(\frac{3}{p_0}-\frac{3}{p}+2)}+1}{C(N)}
( \|U_2\|_{\tilde{L}^\infty_T\dot{B}_{q,1}^{\frac{3}{q}-1}}+\|U_2\|_{L^1_T\dot{B}_{q,1}^{\frac{3}{q}+1}}).
\end{aligned}
\end{equation*}
Collecting the above estimates leads to
\begin{equation}\label{aim2}
\begin{aligned}
\|F_1\|_{L^1_T\dot{B}_{q,1}^{\frac{3}{q}-1}}+\|F_2\|_{L^1_T\dot{B}_{q,1}^{\frac{3}{q}-1}}\le&
\frac{CT^\frac{3}{2}2^{N(\frac{3}{q}+\frac{3}{p_0}-\frac{3}{p}+\frac{5}{2})}}{C(N)^3}
+\frac{CT2^{N(\frac{3}{q}+\frac{3}{2})}}{C(N)}
\|a_1\|_{\tilde{L}^\infty_T\dot{B}_{p_0,1}^\frac{3}{p_0}}\\
&+C
\frac{T2^{N(\frac{3}{p_0}-\frac{3}{p}+2)}+1}{C(N)}
( \|U_2\|_{\tilde{L}^\infty_T\dot{B}_{q,1}^{\frac{3}{q}-1}}+\|U_2\|_{L^1_T\dot{B}_{q,1}^{\frac{3}{q}+1}})\\
&+\|a_1\|_{\tilde{L}^\infty_T\dot{B}_{p_0,1}^\frac{3}{p_0}}\|U_2\|_{L^1_T\dot{B}_{q,1}^{\frac{3}{q}+1}}.
\end{aligned}
\end{equation}
{\bf Step 4. Some estimates of the pressure $P$ and the modified pressure $\Pi$}
By using (\ref{P}),  (\ref{aim1}), (\ref{aim2}) and
$$\|a_0\Delta U_0\|_{L^1_T\dot{B}_{q,1}^{\frac{3}{q}-1}}\le C\|a_0\|_{\dot{B}_{p_0,1}^\frac{3}{p_0}}\|U_0\|_{L^1_T\dot{B}_{q,1}^{\frac{3}{q}+1}}
\le \frac{CT2^{N(\frac{3}{q}+\frac{3}{p_0}-\frac{3}{p}+\frac{3}{2})}}{C(N)^2},
$$
we obtain
\begin{equation}\label{3.4}
\begin{aligned}
\|\nabla P\|_{L^1_T\dot{B}_{q,1}^{\frac{3}{q}-1}}
\le& C\|a\|_{\tilde{L}^\infty_T\dot{B}_{p_0,1}^\frac{3}{p_0}}\|\nabla P\|_{L^1_T\dot{B}_{q,1}^{\frac{3}{q}-1}}
+C\|u\cdot\nabla u\|_{L^1_T\dot{B}_{q,1}^{\frac{3}{q}-1}}+C\|a\Delta u\|_{L^1_T\dot{B}_{q,1}^{\frac{3}{q}-1}}\\
\le& C(\frac{2^{N(\frac{3}{p_0}-\frac{3}{p})}}{C(N)}+\|a_1\|_{\tilde{L}^\infty_T\dot{B}_{p_0,1}^\frac{3}{p_0}})
\|\nabla P\|_{L^1_T\dot{B}_{q,1}^{\frac{3}{q}-1}}+\frac{CT2^{N(\frac{3}{q}+\frac{3}{p_0}-\frac{3}{p}+\frac{3}{2})}}{C(N)^2}\\
&+\frac{CT2^{N(\frac{3}{q}+\frac{3}{2})}}{C(N)}\|a_1\|_{\tilde{L}^\infty_T\dot{B}_{p_0,1}^\frac{3}{p_0}}
+C\frac{1+T2^{N(\frac{3}{p_0}-\frac{3}{p}+2)}}{C(N)}(\|U_2\|_{\tilde{L}^\infty_T\dot{B}_{q,1}^{\frac{3}{q}-1}}\\
&+\|U_2\|_{L^1_T\dot{B}_{q,1}^{\frac{3}{q}+1}})+C\|U_2\|_{\tilde{L}^2_T\dot{B}_{q,1}^\frac{3}{q}}^2
+C\|a_1\|_{\tilde{L}^\infty_T\dot{B}_{p_0,1}^\frac{3}{p_0}}\|U_2\|_{L^1_T\dot{B}_{q,1}^{\frac{3}{q}+1}}.
\end{aligned}
\end{equation}
Similarly, using (\ref{MP}) yields that
\begin{equation}\label{3.5}
\begin{aligned}
\|\nabla \Pi\|_{L^1_T\dot{B}_{q,1}^{\frac{3}{q}-1}}
\le& C\|a\|_{\tilde{L}^\infty_T\dot{B}_{p_0,1}^\frac{3}{p_0}}\|\nabla \Pi\|_{L^1_T\dot{B}_{q,1}^{\frac{3}{q}-1}}\\
&+C\|a\|_{\tilde{L}^\infty_T\dot{B}_{p_0,1}^\frac{3}{p_0}}^2\| u\|_{L^1_T\dot{B}_{q,1}^{\frac{3}{q}+1}}+C\|u\cdot\nabla u\|_{L^1_T\dot{B}_{q,1}^{\frac{3}{q}-1}}\\
\le& C(\frac{2^{N(\frac{3}{p_0}-\frac{3}{p})}}{C(N)}+\|a_1\|_{\tilde{L}^\infty_T\dot{B}_{p_0,1}^\frac{3}{p_0}})\|\nabla \Pi\|_{L^1_T\dot{B}_{q,1}^{\frac{3}{q}-1}}+C\|U_2\|_{\tilde{L}^2_T\dot{B}_{q,1}^\frac{3}{q}}^2\\
&+C(\frac{2^{N(\frac{6}{p_0}-\frac{6}{p})}}{C(N)^2}
+\|a_1\|_{\tilde{L}^\infty_T\dot{B}_{p_0,1}^\frac{3}{p_0}}^2)(\frac{T2^{N(\frac{3}{q}+\frac{3}{2})}}{C(N)}+\|U_2\|_{L^1_T\dot{B}_{q,1}^{\frac{3}{q}+1}})\\
&+\frac{CT2^{N(\frac{3}{q}+1)}}{C(N)^2}+\frac{C}{C(N)}
(\|U_2\|_{\tilde{L}^\infty_T\dot{B}_{q,1}^{\frac{3}{q}-1}}+\|U_2\|_{\tilde{L}^2_T\dot{B}_{q,1}^\frac{3}{q}})
.
\end{aligned}
\end{equation}
{\bf Step 5. The estimate of $\|a_1\|_{\tilde{L}^\infty_T\dot{B}_{p_0,1}^\frac{3}{p_0}}$}
Applying Proposition  \ref{p1.3} to  the transport equation (\ref{TE}), one deduces
$$\|a_1\|_{\tilde{L}^\infty_T\dot{B}_{p_0,1}^\frac{3}{p_0}}\le \|u\cdot\nabla a_0\|_{L^1_T\dot{B}_{p_0,1}^\frac{3}{p_0}}\exp\{C\|\nabla u\|_{L^1_T\dot{B}_{p_0,1}^\frac{3}{p_0}}\}.$$
It follows from using (\ref{kp}) that
$$\|U_0\cdot\nabla a_0\|_{L^1_T\dot{B}_{p_0,1}^\frac{3}{p_0}}\le C(\|U_0\|_{L^1_TL^\infty}\|\nabla a_0\|_{\dot{B}_{p_0,1}^\frac{3}{p_0}}
+\|U_0\|_{L^1_T\dot{B}_{p_0,1}^\frac{3}{p_0}}\|\nabla a_0\|_{L^\infty})\le \frac{CT2^{N(\frac{3}{p_0}-\frac{3}{p}+\frac{3}{2})}}{C(N)^2}
,$$
$$
 \|U_1\cdot\nabla a_0\|_{L^1_T\dot{B}_{p_0,1}^\frac{3}{p_0}}\le C(\|U_1\|_{L^1_TL^\infty}\|\nabla a_0\|_{\dot{B}_{p_0,1}^\frac{3}{p_0}}
+\|U_1\|_{L^1_T\dot{B}_{p_0,1}^\frac{3}{p_0}}\|\nabla a_0\|_{L^\infty})\le \frac{CT^22^{N(\frac{3}{p_0}-\frac{3}{p}+\frac{7}{2})}}{C(N)^3}
.
$$
Using (\ref{p1}), we have
\begin{equation*}
\begin{aligned}
\|U_2\cdot\nabla a_0\|_{L^1_T\dot{B}_{p_0,1}^\frac{3}{p_0}}\le& C(T^\frac{1}{2}\|\nabla a_0\|_{L^\infty}\|U_2\|_{\tilde{L}^2_T\dot{B}_{p_0,1}^\frac{3}{p_0}}
+T\|\nabla^2 a_0\|_{L^\infty}\|U_2\|_{\tilde{L}^\infty_T\dot{B}_{p_0,1}^{\frac{3}{p_0}-1}})\\
\le& \frac{CT^\frac{1}{2}2^{N(1-\frac{3}{p})}}{C(N)}
(\|U_2\|_{\tilde{L}^\infty_T\dot{B}_{p_0,1}^{\frac{3}{p_0}-1}}+\|U_2\|_{\tilde{L}^2_T\dot{B}_{p_0,1}^\frac{3}{p_0}}).
\end{aligned}
\end{equation*}
Using $\dot{B}_{q,1}^{s+\frac{3}{q}}\hookrightarrow \dot{B}_{p_0,1}^{s+\frac{3}{p_0}}$, thus we get
\begin{equation}\label{density}
\begin{aligned}
\|a_1\|_{\tilde{L}^\infty_T\dot{B}_{p_0,1}^\frac{3}{p_0}}\le& \{\frac{CT2^{N(\frac{3}{p_0}-\frac{3}{p}+\frac{3}{2})}}{C(N)^2}+
\frac{CT^\frac{1}{2}2^{N(1-\frac{3}{p})}}{C(N)}
(\|U_2\|_{\tilde{L}^\infty_T\dot{B}_{q,1}^{\frac{3}{q}-1}}+\|U_2\|_{\tilde{L}^2_T\dot{B}_{q,1}^\frac{3}{q}})
\}\\
&\exp\{\frac{CT2^{N(\frac{3}{p_0}+\frac{3}{2})}}{C(N)}+C\|U_2\|_{L^1_T\dot{B}_{p_0,1}^{\frac{3}{p_0}+1}}\}.
\end{aligned}
\end{equation}

\subsection {Proof of Theorem \ref{t1}}
\label{s3.5} Denote
$$X_T:=\|a_1\|_{\tilde{L}^\infty_T\dot{B}_{p_0,1}^\frac{3}{p_0}},\ Y_T:=\|U_2\|_{\tilde{L}^\infty_T\dot{B}_{q,1}^{\frac{3}{q}-1}}+\|U_2\|_{L^1_T\dot{B}_{q,1}^{\frac{3}{q}+1}},$$
and
$$\bar{T}=\sup\left\{ t\in (0,T):\ Y_T\le M_1(\frac{T2^{N(\frac{3}{q}+\frac{6}{p_0}-\frac{6}{p}+\frac{3}{2})}}{C(N)^3}+\frac{1}{C(N)^2}),\ X_T\le M_2\frac{2^{N(\frac{6}{p_0}-\frac{6}{p})}}{C(N)^2}\right\},$$
where $M_i$ $i=1,2$ are large enough  constants, which will be determined later on. Assume $\bar{T}<T$. Choosing $N$ such that
$$\frac{2^{N(\frac{3}{p_0}-\frac{3}{p})}}{C(N)}\ll 1,\ M_2\frac{2^{N(\frac{6}{p_0}-\frac{6}{p})}}{C(N)^2}\ll 1,$$
thanks to (\ref{3.5}) and (\ref{3.4}), we get
\begin{equation*}
\begin{aligned}
\|\nabla \Pi\|_{L^1_T\dot{B}_{q,1}^{\frac{3}{q}-1}}\le& C(\frac{T2^{N(\frac{3}{q}+\frac{6}{p_0}-\frac{6}{p}+\frac{3}{2})}}{C(N)^3}
+\frac{1}{C(N)^2})+\frac{CT2^{N(\frac{3}{q}+\frac{3}{2})}}{C(N)}X_T^2\\
&+\frac{C2^{2N(\frac{3}{p_0}-\frac{3}{p})}}{C(N)^2}Y_T+CX_T^2Y_T
+CY_T^2,
\end{aligned}
\end{equation*}
and
\begin{equation}\label{1003}
\begin{aligned}
\|\nabla P\|_{L^1_T\dot{B}_{q,1}^{\frac{3}{q}-1}}\le&
\frac{CT2^{N(\frac{3}{q}+\frac{3}{p_0}-\frac{3}{p}+\frac{3}{2})}}{C(N)^2}
+\frac{CT2^{N(\frac{3}{q}+\frac{3}{2})}}{C(N)}X_T\\
&+C\frac{T2^{N(\frac{3}{p_0}-\frac{3}{p}+2)}+1}{C(N)}Y_T+CY_T(X_T+Y_T),
\end{aligned}
\end{equation}
respectively.  The estimate (\ref{1003}) yields
\begin{equation*}
\begin{aligned}
\|a\nabla P\|_{L^1_T\dot{B}_{q,1}^{\frac{3}{q}-1}}\le& C(\frac{2^{N(\frac{3}{p_0}-\frac{3}{p})}}{C(N)}+X_T)\|\nabla P\|_{L^1_T\dot{B}_{q,1}^{\frac{3}{q}-1}}\\
\le& \frac{T2^{N(\frac{3}{q}+\frac{6}{p_0}-\frac{6}{p}+\frac{3}{2})}}{C(N)^3}+\frac{CT2^{N(\frac{3}{q}+\frac{3}{p_0}-\frac{3}{p}+\frac{3}{2})}}{C(N)^2}X_T\\
&+C\frac{2^{N(\frac{3}{p_0}-\frac{3}{p})}+T2^{N(\frac{6}{p_0}-\frac{6}{p}+2)}}{C(N)^2}Y_T+\frac{C2^{N(\frac{3}{p_0}-\frac{3}{p})}}{C(N)}Y_T(X_T+Y_T)\\
&+\frac{CT2^{N(\frac{3}{q}+\frac{3}{2})}}{C(N)}X_T^2+CX_TY_T(X_T+Y_T).
\end{aligned}
\end{equation*}
Setting $N$ such that
$$M_1\frac{T^\frac{3}{2}2^{N(\frac{3}{ q}-\frac{3}{p}+\frac{5}{2})}}{C(N)^2}\ll 1,$$
and using (\ref{density}), we have
\begin{equation*}
\begin{aligned}
X_T\le \frac{CT2^{N(\frac{3}{p_0}-\frac{3}{p}+\frac{3}{2})}}{C(N)^2}+\frac{CM_1T^\frac{3}{2}2^{N(\frac{3}{q}+\frac{6}{p_0}-\frac{9}{p}+\frac{5}{2})}}{C(N)^4}
+\frac{CM_1T^\frac{1}{2}2^{N(1-\frac{3}{p})}}{C(N)^3}\le \frac{C_12^{N(\frac{6}{p_0}-\frac{6}{p})}}{C(N)^2}.
\end{aligned}
\end{equation*}
Thanks to the above estimates, choosing $N$ such that
$$\frac{T2^{N(\frac{3}{q}+\frac{6}{p_0}-\frac{6}{p}+\frac{3}{2})}}{C(N)^3}\ll 1,$$
we can obtain
\begin{equation*}
\begin{aligned}
Y_T\le& C\|u\cdot\nabla u\|_{L^1_T\dot{B}_{q,1}^{\frac{3}{q}-1}}+C\sum_{i=1,2}\|F_i\|_{L^1_T\dot{B}_{q,1}^{\frac{3}{q}-1}}+
C\|a\nabla P\|_{L^1_T\dot{B}_{q,1}^{\frac{3}{q}-1}}+C\|\nabla \Pi\|_{L^1_T\dot{B}_{q,1}^{\frac{3}{q}-1}}\\
\le& C(\frac{T2^{N(\frac{3}{q}+\frac{6}{p_0}-\frac{6}{p}+\frac{3}{2})}}{C(N)^3}
+\frac{1}{C(N)^2})+\frac{CT2^{N(\frac{3}{q}+\frac{3}{2})}}{C(N)}X_T\\
&+C(\frac{1+T2^{N(\frac{3}{p_0}-\frac{3}{p}+2)}}{C(N)}+\frac{2^{N(\frac{6}{p_0}-\frac{6}{p})}}{C(N)^2})Y_T
+CY_T(X_T+Y_T)\\
\le& C_2(\frac{T2^{N(\frac{3}{q}+\frac{6}{p_0}-\frac{6}{p}+\frac{3}{2})}}{C(N)^3}
+\frac{1}{C(N)^2})+\frac{1}{2}Y_T,
\end{aligned}
\end{equation*}
which follows
$$Y_T\le 2C_2(\frac{T2^{N(\frac{3}{q}+\frac{6}{p_0}-\frac{6}{p}+\frac{3}{2})}}{C(N)^3}
+\frac{1}{C(N)^2}).$$
We can see from the Remark \ref{rr1} that the conditions in Lemma \ref{l1} can ensure the above requirements of $N$. If we set $M_2=4C_1$ and $M_1=4C_2$, then a contradiction is obtained. Therefore, we have $\bar{T}=T_0$, and
\begin{equation}\label{sm}
Y_T\le 4C_2(\frac{T2^{N(\frac{3}{q}+\frac{6}{p_0}-\frac{6}{p}+\frac{3}{2})}}{C(N)^3}
+\frac{1}{C(N)^2}),\ X_t\le 4\frac{C_12^{N(\frac{6}{p_0}-\frac{6}{p})}}{C(N)^2},\  \forall\ T\le T_0.
\end{equation}
Combining with (\ref{LB}) and (\ref{sm}), we get
$$\|u(T_0)\|_{\dot{B}_{6,1}^{-\frac{1}{2}}}\ge \|U_1(T_0)\|_{\dot{B}_{6,1}^{-\frac{1}{2}}}-(\|U_0(T_0)\|_{\dot{B}_{6,1}^{-\frac{1}{2}}}+Y_{T_0})\ge \frac{c}{8}2^{2(\epsilon_1-\epsilon)N}.$$
Thanks to (\ref{3.1}) and (\ref{100}), then we conclude the proof of Theorem \ref{t1}.
\vskip .4in
\section{Proof of Theorem \ref{t2}}
\label{main2}
The proof  is very similar to the proof of Theorem \ref{t1}.  Let us keep the process in Section \ref{s3.1} in mind. Now we begin with the choice of initial data.
\subsection{The choice of initial data}
 \label{s4.1}
 As Section \ref{s3.2}, we assume (\ref{f1}) and make the same assumptions of $C(N)$, $\epsilon_1$, $\hat{\phi}$. Let us construct the initial data as follows:
\begin{equation*}
\begin{aligned}
\widehat{a_0}(\xi)=&\frac{1}{C(N)}\sum_{k=100}^N2^{-\frac{k}{2}}\left(\hat{\phi}(\xi-2^ke_1)+\hat{\phi}(\xi+2^ke_1)\right),\\
\widehat{u_0}(\xi)=&\frac{1}{C(N)}\sum_{k=100}^N2^{k(1-\frac{3}{p})}\left(\hat{\phi}(\xi+2^ke_1)-\hat{\phi}(\xi-2^ke_1)\right)\frac{1}{|\xi|}
\left(
  \begin{array}{c}
   \xi_2 \\
   -\xi_1 \\
    0 \\
  \end{array}
\right).
\end{aligned}
\end{equation*}
 One can see  $a_0$ and $u_0$ are  real valued function and  real vector-valued function, respectively.  One can also check the following estimates hold, i.e.,
\begin{equation}\label{3.1}
\|a_0\|_{L^\infty}\le \frac{C}{C(N)},\  \|a_0\|_{\dot{B}_{6,1}^\frac{1}{2}}\le \frac{CN}{C(N)},\ \|u_0\|_{\dot{B}_{p,1}^{\frac{3}{p}-1}}\le \frac{CN}{C(N)}
\end{equation}
and
\begin{equation*}
\|a_0\|_{\dot{B}_{r,1}^s}\le \frac{C2^{N(s-\frac{1}{2})}}{C(N)},\ \|u_0\|_{\dot{B}_{r_1,1}^{s_1}}\le \frac{C2^{N(s_1+1-\frac{3}{p})}}{C(N)},
\end{equation*}
with $(r,r_1)\in [1,\infty]^2$ and $s >\frac{1}{2}$, $s_1>\frac{3}{p}-1$.
\vskip .1in
Now, we give a lemma.
\begin{lemma}\label{l2}
Let $p\in(6,\infty)$. Then there exist some positive constants $q$,  $\epsilon$, $\epsilon_1$ satisfying
\begin{equation}\label{cond2}
\left\{
\begin{aligned}
& \frac{3}{q}<\min\{1+2\epsilon-2\epsilon_1,\ \frac{3}{4}+\frac{9}{2p}+2\epsilon-\epsilon_1,\\
&\hspace{15mm} \frac{3}{4}-\frac{3}{2p}+2\epsilon-3\epsilon_1,\
\frac{1}{2}+\frac{3}{p}+2\epsilon-2\epsilon_1\}\\
& q\in (3,6),\ 0<2\epsilon<2\epsilon_1<\frac{1}{2}-\frac{3}{p}.
\end{aligned}
\right.
\end{equation}
\end{lemma}
\begin{rem}\label{rr2}
Although one can easily find that  $1+2\epsilon-2\epsilon_1$ and $\frac{3}{4}+\frac{9}{2p}+2\epsilon-\epsilon_1$ in (\ref{cond2}) can be dropped, here we keep it in the bracket in order to giving a detailed
analysis of the conditions in (\ref{cond20}). Now,
let us  give the following explanations of the limitations in (\ref{cond2}). Let  $T=2^{-2(1+\epsilon)N}$, then we have
\begin{equation}\label{cond20}
\left\{
\begin{aligned}
& \frac{T2^{N(\frac{3}{q}-\frac{3}{p}+\frac{3}{2})}}{C(N)^2}\ll 1\ \ \Longleftrightarrow\ \ \frac{3}{q}<1+2\epsilon-2\epsilon_1,\\
& \frac{T2^{N(\frac{3}{q}-\frac{6}{p}+\frac{3}{2})}}{C(N)}\ll1\ \ \Longleftrightarrow\ \ \frac{3}{q}<\frac{3}{4}+\frac{9}{2p}+2\epsilon-\epsilon_1,\\
& \frac{T2^{N(\frac{3}{q}-\frac{3}{p}+2)}}{C(N)^3}\ll1\ \ \Longleftrightarrow\ \ \frac{3}{q}<\frac{3}{4}-\frac{3}{2p}+2\epsilon-3\epsilon_1,\\
& \frac{T2^{N(\frac{3}{q}-\frac{6}{p}+2)}}{C(N)^2}\ll1\ \ \Longleftrightarrow\ \ \frac{3}{q}<\frac{1}{2}+\frac{3}{p}+2\epsilon-2\epsilon_1,\\
\end{aligned}
\right.
\end{equation}
Actually,  we  assume the conditions on the right hand side of (\ref{cond20}) to  ensure the conditions on the left hand side which is required in  our proof.
Furthermore, we  use $q\in (3,6)$  to ensure some product estimates like (\ref{p3}).
The choice of $C(N)$ needs $2\epsilon_1<\frac{1}{2}-\frac{3}{p}$, while $\epsilon<\epsilon_1$ ensures the norm inflation of $U_1$ (see  Subsection \ref{s4.2}).
\end{rem}
Lemma \ref{l2} can be proved easily, here we use the following  example in this article.
\vskip .1in
If $p\in (6,18]$, the first inequality in (\ref{cond2}) reduces to
$\frac{3}{q}<\frac{3}{4}-\frac{3}{2p}+2\epsilon-3\epsilon_1$, then we can choose
$$\epsilon_1=\frac{1}{4}(\frac{1}{2}-\frac{3}{p}),\ \epsilon=\frac{1}{5}(\frac{1}{2}-\frac{3}{p})$$
and
$$\forall\ \frac{3}{q}\in (\frac{1}{2},\ \frac{23}{40}-\frac{9}{20p}),$$
leading
\begin{equation}\label{300}
2(\epsilon_1-\epsilon)=\frac{1}{10}(\frac{1}{2}-\frac{3}{p}),\ \frac{1}{2}(\frac{1}{2}-\frac{3}{p}-2\epsilon_1)=\frac{1}{4}(\frac{1}{2}-\frac{3}{p});
\end{equation}
\vskip .1in
If $p\in (18,\infty)$, the first inequality in (\ref{cond2}) can be ensured by  $\frac{3}{q}<\frac{1}{2}+\frac{3}{p}+2\epsilon-3\epsilon_1$, then we choose
$$\epsilon_1=\frac{2.1}{p},\ \epsilon=\frac{2}{p},\ \forall\ \frac{3}{q}\in(\frac{1}{2},\frac{1}{2}+\frac{0.7}{p})$$
leading
\begin{equation}\label{30000}
2(\epsilon_1-\epsilon)=\frac{0.2}{p},\ \frac{1}{2}(\frac{1}{2}-\frac{3}{p}-2\epsilon_1)=\frac{1}{4}-\frac{3.6}{p}\ge \frac{0.2}{p}.
\end{equation}
\subsection{The lower bound of $U_1$}
 \label{s4.2}
 Following the same idea and process as Section \ref{s3.3}, one can easily get  the finial lower bound of $\|U_1(t)\|_{\dot{B}_{p,1}^{\frac{3}{p}-1}}$, that is,
\begin{equation}\label{4.1}
\|U_1(t)\|_{\dot{B}_{p,1}^{\frac{3}{p}-1}}\ge  \frac{c}{4}2^{2(\epsilon_1-\epsilon)N},
\end{equation}
where $t=T_0:=2^{-2(1+\epsilon)N}$.
We omit the details to avoid the repetition.
\subsection{The analysis of $U_2$}
 \label{s4.3}
 Let $0\le T\le T_0$, we split this subsection into several steps.\\
{\bf Step 1 The estimate of $U_1$}
Like the previous section, we list first some estimates of $U_1$:
$$\|U_1\|_{L^r_T L^\infty}\le \frac{CT^{1+\frac{1}{r}}2^{N(3-\frac{3}{p})}}{C(N)^2},\
\|U_1\|_{L^1_T\dot{B}_{q_2,1}^{\frac{3}{q_2}+1}}+\|U_1\|_{\tilde{L}^2_T\dot{B}_{q_2,1}^\frac{3}{q_2}}
\le \frac{CT^\frac{3}{2}2^{N(\frac{3}{q_2}-\frac{3}{p}+3)}}{C(N)^2},\ (q_2=6\ {\rm or}\  q).
$$
{\bf Step 2 The estimate of $\|u\cdot\nabla u\|_{L^1_T\dot{B}_{q,1}^{\frac{3}{q}-1}}$}
By using (\ref{kp}),
\begin{equation*}
\begin{aligned}
\|U_0\otimes U_0\|_{L^1_T\dot{B}_{q,1}^\frac{3}{q}}\le& C\|U_0\|_{L^\infty_TL^\infty}\|U_0\|_{L^1_T\dot{B}_{q,1}^\frac{3}{q}}
 \le \frac{CT2^{N(\frac{3}{q}-\frac{6}{p}+2)}}{C(N)^2},\\
\|U_0\otimes U_1\|_{L^1_T\dot{B}_{q,1}^\frac{3}{q}}\le& CT^\frac{1}{2}(\|U_0\|_{L^\infty_TL^\infty}\|U_1\|_{\tilde{L}^2_T\dot{B}_{q,1}^\frac{3}{q}}
+\|U_1\|_{L^\infty_TL^\infty}\|U_0\|_{\tilde{L}^2_T\dot{B}_{q,1}^\frac{3}{q}})\\
 \le& \frac{CT^22^{N(\frac{3}{q}-\frac{6}{p}+4)}}{C(N)^3},\\
\|U_1\otimes U_1\|_{L^1_T\dot{B}_{q,1}^\frac{3}{q}}\le&C\|U_1\|_{L^2_TL^\infty}\|U_1\|_{\tilde{L}^2_T\dot{B}_{q,1}^\frac{3}{q}}\le \frac{CT^32^{N(\frac{3}{q}-\frac{6}{p}+6)}}{C(N)^4},\\
\|U_2\otimes U_2\|_{L^1_T\dot{B}_{q,1}^\frac{3}{q}}\le& C\|U_2\|_{\tilde{L}^2_T\dot{B}_{q,1}^\frac{3}{q}}^2.
\end{aligned}
\end{equation*}
Thanks to (\ref{p1}),
\begin{equation*}
\begin{aligned}
\|U_0\otimes U_2\|_{L^1_T\dot{B}_{q,1}^\frac{3}{q}}\le& C(\|U_0\|_{L^2_TL^\infty}\|U_2\|_{\tilde{L}^2_T\dot{B}_{q,1}^\frac{3}{q}}
+\|\nabla U_0\|_{L^1_TL^\infty}\|U_2\|_{\tilde{L}^\infty_T\dot{B}_{q,1}^{\frac{3}{q}-1}}
)\\
\le& \frac{CT^\frac{1}{2}2^{N(1-\frac{3}{p})}}{C(N)}(\|U_2\|_{\tilde{L}^\infty_T\dot{B}_{q,1}^{\frac{3}{q}-1}}+\|U_2\|_{\tilde{L}^2_T\dot{B}_{q,1}^\frac{3}{q}}).
\end{aligned}
\end{equation*}
\begin{equation*}
\begin{aligned}
\|U_1\otimes U_2\|_{L^1_T\dot{B}_{q,1}^\frac{3}{q}}\le& C(\|U_1\|_{L^2_TL^\infty}\|U_2\|_{\tilde{L}^2_T\dot{B}_{q,1}^\frac{3}{q}}
+\|\nabla U_1\|_{L^1_TL^\infty}\|U_2\|_{\tilde{L}^\infty_T\dot{B}_{q,1}^{\frac{3}{q}-1}}
)\\
\le& \frac{CT^\frac{3}{2}2^{N(3-\frac{3}{p})}}{C(N)^2}(\|U_2\|_{\tilde{L}^\infty_T\dot{B}_{q,1}^{\frac{3}{q}-1}}+\|U_2\|_{\tilde{L}^2_T\dot{B}_{q,1}^\frac{3}{q}}).
\end{aligned}
\end{equation*}
Combining with the above six estimates, we have
\begin{equation}\label{101}
\begin{aligned}
\|u\cdot\nabla u\|_{L^1_T\dot{B}_{q,1}^{\frac{3}{q}-1}}\le&
\frac{CT2^{N(\frac{3}{q}-\frac{6}{p}+2)}}{C(N)^2}
+\frac{CT^\frac{1}{2}2^{N(1-\frac{3}{p})}}{C(N)}\\
&\times(\|U_2\|_{\tilde{L}^\infty_T\dot{B}_{q,1}^{\frac{3}{q}-1}}+\|U_2\|_{\tilde{L}^2_T\dot{B}_{q,1}^\frac{3}{q}})
+C\|U_2\|_{\tilde{L}^2_T\dot{B}_{q,1}^\frac{3}{q}}^2.
\end{aligned}
\end{equation}
{\bf Step 3. The estimate of $\|F_1\|_{L^1_T\dot{B}_{q,1}^{\frac{3}{q}-1}}$ and $\|F_2\|_{L^1_T\dot{B}_{q,1}^{\frac{3}{q}-1}}$}
Using the product estimate (\ref{p3}) and  $\dot{B}_{6,1}^\frac{1}{2}\hookrightarrow \dot{B}_{p_0,1}^\frac{3}{p_0}$,  one gets
\begin{equation*}
\begin{aligned}
\|a_1\Delta U_0\|_{L^1_T\dot{B}_{q,1}^{\frac{3}{q}-1}}\le& C\|a_1\|_{\tilde{L}^\infty_T\dot{B}_{6,1}^\frac{1}{2}}\|U_0\|_{L^1_T\dot{B}_{q,1}^{\frac{3}{q}+1}}\le \frac{CT2^{N(\frac{3}{q}-\frac{3}{p}+2)}}{C(N)}
\|a_1\|_{\tilde{L}^\infty_T\dot{B}_{6,1}^\frac{1}{2}},\\
\|a_1\Delta U_1\|_{L^1_T\dot{B}_{q,1}^{\frac{3}{q}-1}}\le& C\|a_1\|_{\tilde{L}^\infty_T\dot{B}_{6,1}^\frac{1}{2}}\|U_1\|_{L^1_T\dot{B}_{q,1}^{\frac{3}{q}+1}}\le \frac{CT^\frac{3}{2}2^{N(\frac{3}{q}-\frac{3}{p}+3)}}{C(N)^2}
\|a_1\|_{\tilde{L}^\infty_T\dot{B}_{6,1}^\frac{1}{2}},\\
\|a_1\Delta U_2\|_{L^1_T\dot{B}_{q,1}^{\frac{3}{q}-1}}\le& C\|a_1\|_{\tilde{L}^\infty_T\dot{B}_{6,1}^\frac{1}{2}}\|U_2\|_{L^1_T\dot{B}_{q,1}^{\frac{3}{q}+1}},\\
\|a_0\Delta U_1\|_{L^1_T\dot{B}_{q,1}^{\frac{3}{q}-1}}\le& \|a_0\|_{\tilde{L}^\infty_T\dot{B}_{6,1}^\frac{1}{2}}\|U_1\|_{L^1_T\dot{B}_{q,1}^{\frac{3}{q}+1}}
\le\frac{CT^\frac{3}{2}2^{N(\frac{3}{q}-\frac{3}{p}+3)}}{C(N)^3},\\
\|a_0\Delta U_2\|_{L^1_T\dot{B}_{q,1}^{\frac{3}{q}-1}}\le&
 C\|a_0\|_{\dot{B}_{6,1}^\frac{1}{2}}\|U_2\|_{L^1_T\dot{B}_{q,1}^{\frac{3}{q}+1}}
 \le \frac{C}{C(N)}\|U_2\|_{L^1_T\dot{B}_{q,1}^{\frac{3}{q}+1}}.
\end{aligned}
\end{equation*}
Combining  with the above estimates, then we get
\begin{equation}\label{102}
\begin{aligned}
\|F_1\|_{L^1_T\dot{B}_{q,1}^{\frac{3}{q}-1}}+\|F_2\|_{L^1_T\dot{B}_{q,1}^{\frac{3}{q}-1}}\le&
\frac{CT^\frac{3}{2}2^{N(\frac{3}{q}-\frac{3}{p}+3)}}{C(N)^3}+\frac{CT2^{N(\frac{3}{q}-\frac{3}{p}+2)}}{C(N)}
\|a_1\|_{\tilde{L}^\infty_T\dot{B}_{6,1}^\frac{1}{2}}\\
&+\frac{C}{C(N)}\|U_2\|_{L^1_T\dot{B}_{q,1}^{\frac{3}{q}+1}}+
C\|a_1\|_{\tilde{L}^\infty_T\dot{B}_{6,1}^\frac{1}{2}}\|U_2\|_{L^1_T\dot{B}_{q,1}^{\frac{3}{q}+1}}.
\end{aligned}
\end{equation}
{\bf Step 4. The estimates of the pressure $P$ and the modified pressure $\Pi$}
Using the product estimate (\ref{p3}) again, we have
\begin{equation}\label{10000}
\|a_0\Delta U_0\|_{L^1_T\dot{B}_{q,1}^{\frac{3}{q}-1}}\le
C\|a_0\|_{\dot{B}_{6,1}^\frac{1}{2}}\|U_0\|_{L^1_T\dot{B}_{q,1}^{\frac{3}{q}+1}}\le \frac{CT2^{N(\frac{3}{q}-\frac{3}{p}+2)}}{C(N)^2}.
\end{equation}
Using (\ref{P}), (\ref{p3}), (\ref{101}), (\ref{102}) and (\ref{10000}), one obtains
\begin{equation}\label{4.4}
\begin{aligned}
\|\nabla P\|_{L^1_T\dot{B}_{q,1}^{\frac{3}{q}-1}}
\le& C\|a\|_{\tilde{L}^\infty_T\dot{B}_{6,1}^\frac{1}{2}}\|\nabla P\|_{L^1_T\dot{B}_{q,1}^{\frac{3}{q}-1}}
+C\|u\cdot\nabla u\|_{L^1_T\dot{B}_{q,1}^{\frac{3}{q}-1}}+C\|a\Delta u\|_{L^1_T\dot{B}_{q,1}^{\frac{3}{q}-1}}\\
\le& C(\frac{C}{C(N)}+\|a_1\|_{\tilde{L}^\infty_T\dot{B}_{6,1}^\frac{1}{2}})
\|\nabla P\|_{L^1_T\dot{B}_{q,1}^{\frac{3}{q}-1}}\\
&+\frac{CT2^{N(\frac{3}{q}-\frac{3}{p}+2)}}{C(N)^2}+\frac{CT2^{N(\frac{3}{q}-\frac{3}{p}+2)}}{C(N)}\|a_1\|_{\tilde{L}^\infty_T\dot{B}_{6,1}^\frac{1}{2}}\\
&+\frac{CT^\frac{3}{2}2^{N(3-\frac{3}{p})}+C(N)}{C(N)^2}(\|U_2\|_{\tilde{L}^\infty_T\dot{B}_{q,1}^{\frac{3}{q}-1}}+\|U_2\|_{L^1_T\dot{B}_{q,1}^{\frac{3}{q}+1}})\\
&+C(\|a_1\|_{\tilde{L}^\infty_T\dot{B}_{6,1}^\frac{1}{2}}
+\|U_2\|_{\tilde{L}^\infty_T\dot{B}_{q,1}^{\frac{3}{q}-1}})\|U_2\|_{L^1_T\dot{B}_{q,1}^{\frac{3}{q}+1}}.
\end{aligned}
\end{equation}
And using (\ref{MP}), (\ref{101}) and  (\ref{102}), we have
\begin{equation}\label{4.5}
\begin{aligned}
\|\nabla \Pi\|_{L^1_T\dot{B}_{q,1}^{\frac{3}{q}-1}}
\le& C\|a\|_{\tilde{L}^\infty_T\dot{B}_{6,1}^\frac{1}{2}}\|\nabla \Pi\|_{L^1_T\dot{B}_{q,1}^{\frac{3}{q}-1}}\\
&+C\|a\|_{\tilde{L}^\infty_T\dot{B}_{6,1}^\frac{1}{2}}^2\| u\|_{L^1_T\dot{B}_{q,1}^{\frac{3}{q}+1}}+C\|u\cdot\nabla u\|_{L^1_T\dot{B}_{q,1}^{\frac{3}{q}-1}}\\
\le& C(\frac{1}{C(N)}+\|a_1\|_{\tilde{L}^\infty_T\dot{B}_{6,1}^\frac{1}{2}})\|\nabla \Pi\|_{L^1_T\dot{B}_{q,1}^{\frac{3}{q}-1}}+C\|U_2\|_{\tilde{L}^2_T\dot{B}_{q,1}^\frac{3}{q}}^2\\
&+C(\frac{1}{C(N)^2}
+\|a_1\|_{\tilde{L}^\infty_T\dot{B}_{6,1}^\frac{1}{2}}^2)(\frac{T2^{N(\frac{3}{q}-\frac{3}{p}+2)}}{C(N)}+\|U_2\|_{L^1_T\dot{B}_{q,1}^{\frac{3}{q}+1}})\\
&+\frac{CT2^{N(\frac{3}{q}-\frac{6}{p}+2)}}{C(N)^2}+\frac{CT^\frac{1}{2}2^{N(1-\frac{3}{p})}}{C(N)^2}
(\|U_2\|_{\tilde{L}^\infty_T\dot{B}_{q,1}^{\frac{3}{q}-1}}+\|U_2\|_{\tilde{L}^2_T\dot{B}_{q,1}^\frac{3}{q}})
.
\end{aligned}
\end{equation}
{\bf Step 5. The estimate of $\|a_1\|_{\tilde{L}^\infty_T\dot{B}_{6,1}^\frac{1}{2}}$}
Thanks to  Proposition \ref{p1.3}, we have
$$\|a_1\|_{\tilde{L}^\infty_T\dot{B}_{6,1}^\frac{1}{2}}\le \|u\cdot\nabla a_0\|_{L^1_T\dot{B}_{6,1}^\frac{1}{2}}\exp\{C\|\nabla u\|_{L^1_T\dot{B}_{6,1}^\frac{1}{2}}\}.$$
And by (\ref{kp}), one gets
$$\|U_0\cdot\nabla a_0\|_{L^1_T\dot{B}_{6,1}^\frac{1}{2}}\le C(\|U_0\|_{L^1_TL^\infty}\|\nabla a_0\|_{\dot{B}_{6,1}^\frac{1}{2}}
+\|U_0\|_{L^1_T\dot{B}_{6,1}^\frac{1}{2}}\|\nabla a_0\|_{L^\infty})\le \frac{CT2^{N(2-\frac{3}{p})}}{C(N)^2}
,$$
$$
 \|U_1\cdot\nabla a_0\|_{L^1_T\dot{B}_{6,1}^\frac{1}{2}}\le C(\|U_1\|_{L^1_TL^\infty}\|\nabla a_0\|_{\dot{B}_{6,1}^\frac{1}{2}}
+\|U_1\|_{\tilde{L}^2_T\dot{B}_{6,1}^\frac{1}{2}}\|\nabla a_0\|_{L^2_TL^\infty})\le \frac{CT^22^{N(4-\frac{3}{p})}}{C(N)^3}
.
$$
Applying (\ref{p1}), we  can obtain
\begin{equation*}
\begin{aligned}
\|U_2\cdot\nabla a_0\|_{L^1_T\dot{B}_{6,1}^\frac{1}{2}}\le& C(T^\frac{1}{2}\|\nabla a_0\|_{L^\infty}\|U_2\|_{\tilde{L}^2_T\dot{B}_{6,1}^\frac{1}{2}}
+T\|\nabla^2 a_0\|_{L^\infty}\|U_2\|_{\tilde{L}^\infty_T\dot{B}_{6,1}^{-\frac{1}{2}}})\\
\le& \frac{CT^\frac{1}{2}2^\frac{N}{2}}{C(N)}
(\|U_2\|_{\tilde{L}^\infty_T\dot{B}_{6,1}^{-\frac{1}{2}}}+\|U_2\|_{\tilde{L}^2_T\dot{B}_{6,1}^\frac{1}{2}}).
\end{aligned}
\end{equation*}
Thus we get
\begin{equation}\label{density2}
\begin{aligned}
\|a_1\|_{\tilde{L}^\infty_T\dot{B}_{6,1}^\frac{1}{2}}\le&
C\left(\frac{T2^{N(2-\frac{3}{p})}}{C(N)^2}+  \frac{T^\frac{1}{2}2^\frac{N}{2}}{C(N)}
(\|U_2\|_{\tilde{L}^\infty_T\dot{B}_{6,1}^{-\frac{1}{2}}}+\|U_2\|_{\tilde{L}^2_T\dot{B}_{6,1}^\frac{1}{2}})\right)\\
&\times\exp\{C\|\nabla U_2\|_{L^1_T\dot{B}_{6,1}^\frac{1}{2}}\}.
\end{aligned}
\end{equation}
\subsection{Proof of Theorem \ref{t2}}
 \label{s4.4}
Denote
$$X_T:=\|a_1\|_{\tilde{L}^\infty_T\dot{B}_{6,1}^\frac{1}{2}},\ Y_T:=\|U_2\|_{\tilde{L}^\infty_T\dot{B}_{q,1}^{\frac{3}{q}-1}}+\|U_2\|_{L^1_T\dot{B}_{q,1}^{\frac{3}{q}+1}},$$
and
$$\bar{T}:=\sup\{t\in (0,T_0):\ Y_T\le M_3T(\frac{2^{N(\frac{3}{q}-\frac{3}{p}+2)}}{C(N)^3}+\frac{2^{N(\frac{3}{q}-\frac{6}{p}+2)}}{C(N)^2}),\ X_T\le M_2\frac{T^\frac{1}{2}2^N}{C(N)^2}\},$$
where $M_3$ and $M_4$ will be fixed later. Using (\ref{density2}), choosing $N$ such that
$$M_3T\left\{\frac{2^{N(\frac{3}{q}-\frac{3}{p}+\frac{3}{2})}}{C(N)^2}+\frac{2^{N(\frac{3}{q}-\frac{6}{p}+\frac{3}{2})}}{C(N)}\right\}\ll 1,$$
we have
\begin{equation}\label{aim3}
\begin{aligned}
X_T\le \frac{CT^\frac{1}{2}2^N}{C(N)^2}+\frac{CT^\frac{1}{2}2^N}{C(N)^2}
\{M_3T(\frac{2^{N(\frac{3}{q}-\frac{3}{p}+\frac{3}{2})}}{C(N)^2}+\frac{2^{N(\frac{3}{q}-\frac{6}{p}+\frac{3}{2})}}{C(N)})\}\le \frac{C_1T^\frac{1}{2}2^N}{C(N)^2}.
\end{aligned}
\end{equation}
From the estimate (\ref{4.4}), one has
\begin{equation*}
\begin{aligned}
\|\nabla P\|_{L^1_T\dot{B}_{q,1}^{\frac{3}{q}-1}}
\le& \frac{1}{2}\|\nabla P\|_{L^1_T\dot{B}_{q,1}^{\frac{3}{q}-1}}+\frac{CT2^{N(\frac{3}{q}-\frac{3}{p}+2)}}{C(N)^2}\\
&+\frac{CT^\frac{3}{2}2^{N(3-\frac{3}{p})}+C(N)}{C(N)^2}Y_T\\
&+\frac{CT2^{N(\frac{3}{q}-\frac{3}{p}+2)}}{C(N)}X_T+CY_T(X_T+Y_T),
\end{aligned}
\end{equation*}
which leads to
\begin{equation*}
\begin{aligned}
\|\nabla P\|_{L^1_T\dot{B}_{q,1}^{\frac{3}{q}-1}}
\le& \frac{CT2^{N(\frac{3}{q}-\frac{3}{p}+2)}}{C(N)^2}
+\frac{CT^\frac{3}{2}2^{N(3-\frac{3}{p})}+C(N)}{C(N)^2}Y_T\\
&+\frac{CT2^{N(\frac{3}{q}-\frac{3}{p}+2)}}{C(N)}X_T+CY_T(X_T+Y_T),
\end{aligned}
\end{equation*}
Applying product estimate again, and using $X_T\ll 1$,  we can deduce that
\begin{equation}\label{5.1}
\begin{aligned}
\|a\nabla P\|_{L^1_T\dot{B}_{q,1}^{\frac{3}{q}-1}}
\le& \frac{CT2^{N(\frac{3}{q}-\frac{3}{p}+2)}}{C(N)^3}+\frac{CT^\frac{3}{2}2^{N(3-\frac{3}{p})}+C(N)}{C(N)^2}Y_T\\
&+\frac{CT2^{N(\frac{3}{q}-\frac{3}{p}+2)}}{C(N)}X_T+CY_T(X_T+Y_T).
\end{aligned}
\end{equation}
Similarly, thanks to (\ref{4.5}), we have
\begin{equation}\label{5.2}
\begin{aligned}
\|\nabla \Pi\|_{L^1_T\dot{B}_{q,1}^{\frac{3}{q}-1}}
\le& CT(\frac{2^{N(\frac{3}{q}-\frac{3}{p}+2)}}{C(N)^3}+\frac{2^{N(\frac{3}{q}-\frac{6}{p}+2)}}{C(N)^2})\\
&+\frac{CT2^{N(\frac{3}{q}-\frac{3}{p}+2)}}{C(N)}X_T^2+(\frac{C}{C(N)^2}+Y_T)Y_T+CX_T^2Y_T.
\end{aligned}
\end{equation}
Collecting the above estimates (\ref{101}), (\ref{102}), (\ref{5.1}), (\ref{5.2}) and (\ref{aim3}), and setting $N$ such that
$$M_3T(\frac{2^{N(\frac{3}{q}-\frac{3}{p}+2)}}{C(N)^3}+\frac{2^{N(\frac{3}{q}-\frac{6}{p}+2)}}{C(N)^2})\ll 1,$$
we have
\begin{equation*}
\begin{aligned}
Y_T\le& C(\|u\cdot\nabla u\|_{L^1_T\dot{B}_{q,1}^{\frac{3}{q}-1}}+\sum_{i=1,2}\| F_i\|_{L^1_T\dot{B}_{q,1}^{\frac{3}{q}-1}}
+\|a\nabla P\|_{L^1_T\dot{B}_{q,1}^{\frac{3}{q}-1}}+\|\nabla \Pi\|_{L^1_T\dot{B}_{q,1}^{\frac{3}{q}-1}}
)\\
\le& C_2T(\frac{2^{N(\frac{3}{q}-\frac{3}{p}+2)}}{C(N)^3}+\frac{2^{N(\frac{3}{q}-\frac{6}{p}+2)}}{C(N)^2})+\frac{1}{2}Y_T.
\end{aligned}
\end{equation*}
This yields
$$Y_T\le 2C_2T(\frac{2^{N(\frac{3}{q}-\frac{3}{p}+2)}}{C(N)^3}+\frac{2^{N(\frac{3}{q}-\frac{6}{p}+2)}}{C(N)^2})$$
One can see from the Remark \ref{rr2} that the conditions in Lemma \ref{l2} can ensure the above requirements of $N$.
Choosing  $M_3=4C_1$ and $M_4=4C_2$, we can get a contradiction by using the continuation argument.
 Therefore, we have $\bar{T}=T_0$, and
\begin{equation}\label{sm2}
Y_T\le 4 C_2T(\frac{2^{N(\frac{3}{q}-\frac{3}{p}+2)}}{C(N)^3}+\frac{2^{N(\frac{3}{q}-\frac{6}{p}+2)}}{C(N)^2}),\ X_t\le 4C_1\frac{T^\frac{1}{2}2^N}{C(N)^2},\  \forall\ T\le T_0.
\end{equation}
Combining with (\ref{4.1}) and (\ref{sm2}), we get
$$\|u(T_0)\|_{\dot{B}_{p,1}^{\frac{3}{p}-1}}\ge \|U_1(T_0)\|_{\dot{B}_{p,1}^{\frac{3}{p}-1}}-(\|U_0(T_0)\|_{\dot{B}_{p,1}^{\frac{3}{p}-1}}+Y_{T_0})\ge \frac{c}{8}2^{2(\epsilon_1-\epsilon)N}.$$
Thanks to (\ref{3.1}), (\ref{300}) and (\ref{30000}), we can complete the proof of Theorem
\ref{t2}.
\vskip .3in
\appendix
\section{The endpoint case for Theorem \ref{t2}: $p=\infty$}
\label{app}
In this section, we give some comments on the endpoint case for Theorem \ref{t2}, that is, the case $p=\infty$. The previous proof in Section \ref{main2}
 is not suit for this case, since some difficulties occur when we bound  $\|U_0\cdot\nabla U_0\|_{L^1_T\dot{B}_{q,1}^{\frac{3}{q}-1}}$ in the estimate of $Y_T$. More precisely, we have
$$\|U_0\cdot\nabla U_0\|_{L^1_T\dot{B}_{q,1}^{\frac{3}{q}-1}}\le \frac{CT2^{N(\frac{3}{q}+2)}}{C(N)^2},$$
combined with (\ref{cond2}) and (\ref{cond20}) yielding
$$\frac{3}{q}<\frac{1}{2}+2\epsilon-2\epsilon_1<\frac{1}{2}.$$
This is a contradiction with $q\in (3,6)$. However, we can brake this barrier by introducing another new modified pressure
$$\Pi_1:=P+\mu|D|^{-2}{\rm div}(a\Delta u)-|D|^{-2}{\rm div}(u\cdot\nabla u)$$
leading
\begin{equation*}
\begin{aligned}
\nabla \Pi_1=&|D|^{-2}\nabla {\rm div}(a\nabla \Pi_1)-\mu|D|^{-2}\nabla {\rm div}\{a|D|^{-2}\nabla {\rm div}(a\Delta u)\}\\
&+\mu|D|^{-2}\nabla {\rm div}\{a|D|^{-2}\nabla {\rm div}(u\cdot\nabla u)\}.
\end{aligned}
\end{equation*}
Then we consider
\begin{equation} \label{INS30}
\left\{
\begin{aligned}
& \partial_t a + u\cdot\nabla a= 0,\   {\rm div}u=0 \\
& \partial_t  u -\mu \Delta u =-u\cdot\nabla u-\nabla \Pi_1-a\nabla P\\
&\hspace{15mm}+\mu a\Delta u+\mu |D|^{-2}\nabla{\rm div}(a\Delta u)-\nabla |D|^{-2}{\rm div}(u\cdot\nabla u),\\
& (a(0,x), u(0,x))=(a_0(x), u_0(x)),
\end{aligned}
\right.
\end{equation}
and have
\begin{equation*}
\begin{aligned}
u(t,x)=&e^{\mu\Delta t}u_0+\int_0^t e^{\mu \Delta(t-\tau)}\{-u\cdot\nabla u-\nabla \Pi_1-a\nabla P\\
&+\mu a\Delta u+\mu |D|^{-2}\nabla{\rm div}(a\Delta u)-\nabla |D|^{-2}{\rm div}(u\cdot\nabla u)\}d\tau.
\end{aligned}
\end{equation*}
Denote
\begin{equation*}
\begin{aligned}
U_0(t)=&e^{\mu \Delta t}u_0,\\
U_1(t)=&\mu\int_0^t e^{\mu \Delta (t-\tau)}\{\underbrace{a_0\Delta U_0+|D|^{-2}\nabla {\rm div}(a_0\Delta U_0)}_{\Xi_1}\underbrace{-U_0\cdot\nabla U_0-\nabla |D|^{-2}{\rm div}(U_0\cdot\nabla U_0)}_{\Upsilon_1} \}d\tau,\\
U_2(t)=&\int_0^t e^{\mu \Delta (t-\tau)}\{K_1+K_2-u\cdot\nabla u-\nabla \Pi_1-a \nabla P+F_1+F_2\}d\tau,
\end{aligned}
\end{equation*}
where
\begin{equation*}
\begin{aligned}
K_1=&(U_1+U_2)\cdot\nabla (U_0+U_1+U_2)+U_0\cdot\nabla (U_1+U_2),\\
K_2=&\nabla |D|^{-2}{\rm div}\{(U_1+U_2)\cdot\nabla (U_0+U_1+U_2)+U_0\cdot\nabla (U_1+U_2)\},\\
F_1=&\mu (a_1\Delta u+a_0\Delta(U_1+U_2)),\ \ F_2=\mu|D|^{-2}\nabla {\rm div}\{a_1\Delta u+a_0\Delta (U_1+U_2)\}.
\end{aligned}
\end{equation*}
Thus, we have the decomposition of $u$. We also  use the previous decomposition of $a$ in the following.   We choose the initial data as the Section \ref{main2} by setting $p=\infty$. The  proof is very  similar, here we only show the framework.
\vskip .1in
Firstly, one can get the large lower bound of $U_1$, which can be obtained  from the estimate of $\mu\int_0^t e^{\mu \Delta (t-\tau)} \Xi_1d\tau$. In fact, the associated bound of $\mu\int_0^t e^{\mu \Delta (t-\tau)}\Upsilon_1 d\tau$ can be absorbed by this lower bound.
\vskip .1in
Secondly, thanks to the modified pressure $\Pi_1$, we can avoid  the estimate of  $\|U_0\cdot\nabla U_0\|_{L^1_T\dot{B}_{q,1}^{\frac{3}{q}-1}}$ and get the small bound of $X_T$ and $Y_T$ by following the procedure as Section \ref{main2}.
\vskip .1in
At last, combining with the above  arguments yields  the desired result by following the same procedure as Section \ref{main2}.
\section{Proof of the existence of the solution}
In this section, we give a brief structure to show the existence of the solution.  In fact, (\ref{INS3}) is equal to
\begin{equation} \label{INS10}
\left\{
\begin{aligned}
& \partial_t a + u\cdot\nabla a= 0,\   {\rm div}u=0 \\
& u(t,x)=e^{\mu\Delta t}u_0+\int_0^t e^{\mu \Delta(t-\tau)}\{-u\cdot\nabla u-\nabla \Pi-a\nabla P\\
&+\mu a\Delta u+\mu |D|^{-2}\nabla{\rm div}(a\Delta u)\}d\tau,\\
& (a(0,x), u(0,x))=(a_0(x), u_0(x)).
\end{aligned}
\right.
\end{equation}
From the previous parts, one can see $U_0$ and $U_1$ are chose as follows:
$$U_0(t)=e^{\mu \Delta t}u_0,$$
$$U_1(t)=\mu\int_0^t e^{\mu \Delta (t-\tau)}\{a_0\Delta U_0+|D|^{-2}\nabla {\rm div}(a_0\Delta U_0) \}d\tau.$$
In addition,  we can get from section \ref{s3.5} that  $(a_1, U_2)$ is the solution to the system below
\begin{equation*} 
\left\{
\begin{aligned}
& \partial_t a_1 + (U_0+U_1+U_2)\cdot\nabla (a_0+a_1)= 0,\   {\rm div}u=0 \\
& U_2(t)=\int_0^t e^{\mu \Delta (t-\tau)}\\
&\ \ \ \ \{-(U_0+U_1+U_2)\cdot\nabla(U_0+U_1+U_2)-\nabla \Pi-(a_0+a_1) \nabla P+F_1+F_2\}d\tau,\\
& (a(0,x), U_2(0,x))=(a_0(x), 0),
\end{aligned}
\right.
\end{equation*}
where
$$F_1=\mu (a_1\Delta u+a_0\Delta(U_1+U_2)),$$
$$F_2=\mu|D|^{-2}\nabla {\rm div}\{a_1\Delta (U_0+U_1+U_2)+a_0\Delta (U_1+U_2)\},$$
and
\begin{equation*}
\begin{aligned}
\nabla \Pi=&|D|^{-2}\nabla {\rm div}((a_0+a_1)\nabla \Pi)\\
&-\mu |D|^{-2}\nabla {\rm div}\left\{(a_0+a_1)|D|^{-2}\nabla {\rm div}((a_0+a_1)\Delta (U_0+U_1+U_2))\right\}\\
&+|D|^{-2}\nabla {\rm div}((U_0+U_1+U_2)\cdot\nabla (U_0+U_1+U_2)),
\end{aligned}
\end{equation*}
although section \ref{s3.5} only gives a priori estimate. The strict proof can follow the Chapter 10 in \cite{BCD}, which is very standard, so we omit the details. Therefore, one can see $a:=a_0+a_1$ and $u:=U_0+U_1+U_2$ is a solution to 
(\ref{INS10}).


\end{document}